\documentclass[
leqno,          
12pt,           
a4paper         
]{amsart}

\usepackage[colorlinks,citecolor=blue,urlcolor=blue]{hyperref}          
\usepackage{courier}                                                    
\usepackage{amssymb, amsmath, amsfonts, amsthm}                         
\usepackage{mathtools, cite, enumerate, rotating, framed, lipsum}       
\usepackage{fancybox, bbm}                                              
\usepackage[dvipsnames]{xcolor}                                         
\usepackage[polish, english]{babel}                                     
\usepackage[T1]{fontenc}                                                
\usepackage[cp1250]{inputenc}                                           
\usepackage[normalem]{ulem}
\usepackage{yfonts}
\usepackage{mathrsfs}
\usepackage{float}

\usepackage[a4paper,top=3cm,bottom=2.5cm,left=3cm,right=3cm,marginparwidth=1.75cm]{geometry}

\usepackage{amsmath,amstext,amssymb,amsopn,amsthm,enumitem}
\usepackage{graphicx}
\usepackage[colorinlistoftodos]{todonotes}
\usepackage{yfonts}
\usepackage{mathrsfs}
\usepackage{xstring}

\DeclareMathAlphabet{\mathpzc}{OT1}{pzc}{m}{it}
 \numberwithin{equation}{section}                        

\newcommand{\thmcount}{equation}                 
\newcounter{specialcounter}

\newtheorem{Thm}[\thmcount]{Theorem}
\newtheorem{Sthm}[specialcounter]{Theorem}
\newtheorem{Cor}[\thmcount]{Corollary}
\newtheorem{Lem}[\thmcount]{Lemma}
\newtheorem{Prop}[\thmcount]{Proposition}
\newtheorem{Rem}[\thmcount]{Remark}
\newtheorem{Defn}[\thmcount]{Definition}
\newtheorem{Ex}[\thmcount]{Example}
\newtheorem{Asu}[\thmcount]{Assumption}
\newtheorem{Sol}[\thmcount]{Solution}
\newtheorem*{Thmx}{Theorem}
\newtheorem*{Corx}{Corollary}
\newtheorem*{Lemx}{Lemma}
\newtheorem*{Propx}{Proposition}
\newtheorem*{Remx}{Remark}
\newtheorem*{Defnx}{Definition}
\newtheorem*{Exx}{Example}
\newtheorem*{Asux}{Assumption}
\newtheorem*{Solx}{Solution}

\newcommand \eq[1]{\begin{equation} #1 \end{equation}}
\newcommand \eqx[1]{\begin{equation*}  #1 \end{equation*}}

\newcommand \alx[1]{\begin{align*}  #1 \end{align*}}
\renewcommand \sp[1]{\begin{equation} \begin{split} #1 \end{split} \end{equation}}
\newcommand \spx[1]{\begin{equation*} \begin{split} #1 \end{split} \end{equation*}}
\newcommand \en[1]{\begin{enumerate}  #1 \end{enumerate}}
\newcommand \ite[1]{\begin{itemize}  #1 \end{itemize}}
\newcommand{\thm}[2]{\begin{Thm} \label{#1} #2 \end{Thm}}

\newcommand{\lem}[2]{\begin{Lem} \label{#1} #2 \end{Lem}}

\newcommand{\cor}[2]{\begin{Cor} \label{#1} #2 \end{Cor}}

\newcommand{\prop}[2]{\begin{Prop} \label{#1} #2 \end{Prop}}

\newcommand{\rem}[2]{\begin{Rem} \label{#1} #2 \end{Rem}}

\newcommand{\defn}[2]{\begin{Defn} \label{#1} #2 \end{Defn}}

\newcommand{\pr}[1]{\begin{proof} #1 \end{proof}}

\newcounter{comcount}

\renewcommand{\hline}{\vbox{\hrule width\textwidth height 1pt}\smallskip}

\renewcommand{\a}{\alpha}                \newcommand{\e}{\varepsilon}
                 
        \newcommand{\la}{\lambda}

\newcommand{\BB}{\mathbb{B}}  
\newcommand{\CC}{\mathbb{C}}  
  \newcommand{\dd}{\mathcal{D}}

\newcommand{\KK}{\mathbb{K}}

  \newcommand{\qq}{\mathcal{Q}}
\newcommand{\RR}{\mathbb{R}}  
  
\newcommand{\TT}{\mathbb{T}} \newcommand{\Tt}{\mathbf{T}}

\newcommand{\WW}{\mathbb{W}}

\newcommand{\ZZ}{\mathbb{Z}}  

\newcommand{\supp}{\mathrm{supp}}
\newcommand{\8}{\infty}

\renewcommand{\rm}[1]{\mathrm{#1}}
\newcommand{\wt}[1]{\widetilde{#1}}

\newcommand{\abs}[1]{\left| #1 \right|}
\newcommand{\set}[1]{\left\{ #1 \right\}}
\newcommand{\norm}[1]{\left\| #1 \right\|}

\newcommand{\eee}[1]{\left( #1 \right)}

\newcommand{\BE}{\WW_{t,\nu}^{\rm{exo}}}
\newcommand{\BC}[1]{\WW_{t,#1}^{\rm{cls}}}
\renewcommand{\wt}{\widetilde}

\newcommand{\tti}[1]{T_t^{[#1]}}
\newcommand{\s}[2]{%
    \IfEqCase{#2}{%
        {1}{#1^*}%
        {2}{#1^{**}}%
        {3}{#1^{***}}
    }
}

\renewcommand{\Tt}{T_t}

\newcommand{\Bcls}{\BB_\nu^{\rm{cls}}}
\newcommand{\Bexo}{\BB_\nu^{\rm{exo}}}

\newcommand{\mx}{\mathbbm{x}}

\newcommand{\my}{\mathbbm{y}}

\DeclareMathOperator*{\esssup}{ess\,sup}

\title{The atomic Hardy space for a general Bessel operator }

\author[ Edyta Kania-Strojec ]{ Edyta Kania-Strojec }
 \address{
 Edyta Kania-Strojec \newline
 \indent Instytut Matematyczny, Uniwersytet Wroc\l awski \newline
 \indent pl. Grunwaldzki 2/4, 50-384 Wroc\l aw, Poland }
  \email{edyta.kania-strojec@math.uni.wroc.pl}

\thanks{ Research is supported by the grant OPUS 2017/25/B/ST1/00599 from the National Science Centre (Narodowe Centrum Nauki), Poland. }

\begin{document}
\maketitle

\begin{abstract}
We study Hardy spaces associated with a general multidimensional Bessel operator $\BB_\nu$. This 
operator depends on a multiparameter of type $\nu$ that is usually restricted to a product of half-lines. Here we deal with the Bessel operator in the general context, with no restrictions on the type parameter. We define the Hardy space $H^1$ for $\BB_\nu$ in terms of the maximal operator of the semigroup of operators $\exp(-t\BB_\nu)$. Then we prove that, in general, $H^1$ admits an atomic decomposition of local type.
\end{abstract}

\section{Introduction} Let $\nu = (\nu_1, \ldots ,\nu_d) \in \RR^d$, $d\geq 1$, and consider the multidimensional Bessel differential operator
\eq{\label{multi-bessel}
\BB_\nu f(x) = \sum_{j=1}^d -\partial_j^2 f(x) - \frac{2\nu_j+1}{x_j} \partial_j f(x),
}
acting on functions on $\RR_+^d=(0,\8)^d$. 
Operator $\BB_\nu$ is formally symmetric in $L^2(\RR_+^d, d\mu_\nu)$, where
\eqx{\label{measure}
 d\mu_\nu(x) = x^{2\nu+1} \, dx :=  x_1^{2\nu_1+1} \ldots x_d^{2\nu_d+1} \, dx_1 \ldots dx_d.}

\subsection{Hardy spaces associated with the classical Bessel operator.} When $\nu \in (-1, \8)^d$ there exists a classical self-adjoint extension of $\BB_\nu$ (acting initially on $C_c^2(\RR_+^d)$), from now on denoted by $\Bcls$ \cite{MR199636, Nowak_Sjogren_Szarek-2017}. The operator $\Bcls$ is the infinitesimal generator of the classical Bessel semigroup of operators $\BC{\nu} = \exp(-t \Bcls)$, which has the integral representation $\BC{\nu} f(x) = \int_{\RR_+^d} \BC{\nu}(x,y) f(y) \, d\mu_\nu(y)$, $t>0$. It is well known that
\eqx{\label{classical-kernel}
\BC{\nu}(x,y) = \prod_{j=1}^d \frac{1}{2t}(x_jy_j)^{-\nu_j} I_{\nu_j}\left(  \frac{x_jy_j}{2t}\right) \exp\left( - \frac{x_j^2+y_j^2}{4t} \right) =: \prod_{j=1}^d W_{t,\nu_j}^{\rm{cls}}(x_j,y_j),
}
for  $x,y\in \RR_+^d$ and $t>0$, where $I_\tau$ denotes the modified Bessel function of the first kind and order $\tau > -1$, cf.\ \cite{Watson}. The classical Bessel kernel satisfies the lower and upper Gaussian bounds, i.e.
\eq{\label{classical-gauss}
\frac{C_1}{\mu_\nu(B(x,\sqrt{t}))} \exp\left( -\frac{|x-y|^2}{c_1t}\right) \leq \BC{\nu}(x,y) \leq \frac{C_2}{\mu_\nu(B(x,\sqrt{t}))} \exp\left( -\frac{|x-y|^2}{c_2t}\right),
}
with some constants $c_1,c_2, C_1, C_2 >0$, where $B(x,\sqrt{t}) = \set{y\in\RR_+^d \ : \ |x-y| < \sqrt{t}}$.

Recently, harmonic analysis related to the classical Bessel operator has been extensively developed, see e.g.\ \cite{BCN,BDT_d'Analyse,BCC,Nowak_Sjogren_Szarek-2017,EK_MP_JFAA_2019,DPW_JFAA,MR2345777, MR3780134,MR3176917,MR3631457,MR2600427} and references therein. In particular, the Hardy space
\eq{\label{h1-classical}
H^1(\Bcls) = \set{f\in L^1(\RR_+^d, d\mu_\nu) \ : \ \norm{f}_{H^1(\Bcls)} : = \norm{\sup_{t>0} \abs{\BC{\nu}f}}_{L^1(\RR_+^d, d\mu_\nu)} < \8},}
associated with $\Bcls$ and its characterizations have been studied. 
An especially useful result is the characterization of $H^1(\Bcls)$ by atomic decomposition, which was proved in \cite{BDT_d'Analyse} in the one-dimensional case and then extended to higher dimensions in \cite{DPW_JFAA}. This atomic characterization also follows from a more general result from \cite{Dziubanski2017}. This result states, in particular, that every function $f\in H^1(\Bcls)$ can be represented as $f = \sum_k \la_k a_k$, where $\sum_k \abs{\la_k}\simeq \norm{f}_{H^1(\Bcls)}$ and $a_k$ are the classical atoms on the space of homogeneous type $(\RR_+^d, |\cdot|, d\mu_\nu)$, that is there exist balls $B_k$ such that
\eq{ \label{classical-atoms}
\supp \ a_k \subseteq B_k, \quad \norm{a_k}_\8 \leq \mu_\nu(B_k)^{-1}, \quad \int a_k \,d\mu_\nu = 0.
}
We say that the atoms $a_k$ satisfy localization, size and cancellation conditions \cite{CoifmanWeiss_BullAMS}.

The main goal of this paper is to prove an atomic decomposition theorem for the Hardy space associated with the multidimensional Bessel operator $\BB_\nu$  in the general situation, admitting all $\nu\in\RR^d$. For a precise definition of the operator $\BB_\nu$ see in Section \ref{secion-main-results}. According to the best of the author's knowledge the theory of Hardy spaces for the Bessel operator for the full range of the parameter has never been studied. The most general results are known for $\nu_j > -1$, and a vast majority of them are restricted to $\nu_j\geq -1/2$, $j=1,\ldots,d$. Our atomic decomposition theorem is a starting point for developing this theory in a general Bessel context and the first step towards proving other characterizations of $H^1(\BB_\nu)$. Note that one of the difficulties in the general situation is that the measure is not locally finite, in particular it does not satisfy the standard doubling condition. Theory of Hardy spaces or Calder\'on-Zygmund operators for spaces with non-doubling measure is far more difficult and less known. It has been developed only quite recently, see e.g. \cite{MR1812821, MR1928090, MR3157341,MR1626935,MR1998349}. Usually an essential assumption in this theory is the polynomial growth condition $\mu(B(x,r)) \leq C r^n$. But this condition fails for $\mu_\nu$ unless $\nu\in(-1,\8)^d$, since otherwise there are balls of arbitrarily small radius and infinite measure. For the same reason, our measure does not satisfy the property known in the literature as a local doubling condition, considered for instance in \cite{MR2581426}. 
 Another obstacle is that, in general, operators considered in the literature are generators of semigroups of operators whose integral kernels satisfy the upper Gaussian bounds. This is not the case of $\BB_\nu$ if $\nu\notin (-1,\8)^d$. Nevertheless we overcome these difficulties and give an atomic decomposition of the Hardy space associated with the general multidimensional Bessel operator $\BB_\nu$ for any $\nu\in\RR^d$.

\subsection{The exotic Bessel operator}
Recall that $d\mu(x) = x^{2\nu+1} dx$. When we take $\nu\in((-\8,1)\setminus\set{0})^d$ in \eqref{multi-bessel}, then there exists a self-adjoint extension of $\BB_\nu$, denoted by $\Bexo$ and called the exotic Bessel operator \cite{Nowak_Sjogren_Szarek-2017}. It is quite remarkable that for the parameter $\nu\in((-1,1)\setminus\set{0})^d$ there exist both, the classical and the exotic, self-adjoint extensions of $\BB_\nu$, which differ significantly. We denote by $\BE = \exp(-t\Bexo)$ the semigroup of operators generated by $\Bexo$. It is known that the exotic semigroup kernel can be expressed in a simple way in terms of the classical one, namely
\spx{
	\BE(x,y) = \prod_{j=1}^d W_{t,\nu_j}^{\rm{exo}}(x_j,y_j) = \prod_{j=1}^d (x_jy_j)^{-2\nu_j}W_{t,-\nu_j}^{\rm{cls}}(x_j,y_j).
}
It turns out that if $0<\nu_j<1$ then the so-called pencil phenomenon occurs. In the one-dimensional situation this means that for any fixed $t>0$ the operator $W_{t,\nu_j}^{\rm{exo}}$ is well defined on $L^p(\RR_+,d\mu_{\nu_j})$ and maps this space into itself if and only if $\nu_j+1<p<(\nu_j+1)/\nu_j$. Therefore there is no reason in studying 
the theory of Hardy spaces for the exotic Bessel operator in the case $\nu\in(0,1)^d$. For more details see the discussion in \cite[Sec.\ 4]{Nowak_Sjogren_Szarek-2017} or in \cite{MR3780134}. Thus, from now on we shall consider the operator $\Bexo$ only for $\nu\in(-\8,0)^d$.

For the sake of convenience we shall write $\BB_{-\nu}^{\rm{exo}}$, where $\nu\in(0,\8)^d$. Hence, we will use such notation in the rest of the paper.


\subsection{Main results} \label{secion-main-results}
We consider a general multidimensional Bessel operator $\BB_\nu=L_1 + \ldots + L_d$, $d\geq 1$, where each $L_i$ is either one-dimensional classical Bessel operator $B_{\nu_i}^{\rm{cls}}$ for $\nu_i\in(-1,\8)$ or one-dimensional exotic Bessel operator $B_{-\nu_i}^{\rm{exo}}$ for $\nu_i \in(0,\8)$ (acting on the $i$th coordinate variable, see Section \ref{sec-product-1} for details). Since $B_{\nu_i}^{\rm{cls}}$ and $B_{-\nu_i}^{\rm{exo}}$ are self-adjoint operators, $\BB_\nu$ is well defined and essentially self-adjoint \cite[Thm.\ 7.23]{MR2953553}. To simplify the notation, by changing the coordinates, we shall consider $\BB_\nu = \BB_{\nu_c}^{\rm{cls}} + \BB_{-\nu_e}^{\rm{exo}}$, where $\nu=(\nu_c,-\nu_e)\in(-1,\8)^{d_1}\times (-\8,0)^{d_2}$ and $d_1+d_2=d$. Then
\eqx{
	\BB_\nu f(\mx)= \BB_{\nu_c}^{\rm{cls}} f(\cdot,\mx_2) + \BB_{-\nu_e}^{\rm{exo}}f(\mx_1,\cdot) \qquad \text{for} \quad \mx = (\mx_1,\mx_2) \in \RR_+^{d_1} \times \RR_+^{d_2} = \RR_+^d.
}
As we mentioned above, the case $d_1\geq1$ and $d_2=0$, that is the classical case, is well known and slightly different from the context involving the exotic Bessel operator. Therefore we will consider $d_1\geq 0$ and $d_2\geq 1$. 

Denote by $\WW_{t,\nu} = \exp(-t\BB_\nu)$ the semigroup of operators generated by $\BB_\nu$. Clearly, the semigroup $\WW_{t,\nu}$ has the integral representation 
\eqx{
	\WW_{t,\nu} f(\mx) = \int_{\RR_+^d} \WW_{t,\nu}(\mx,\my) f(\my) \,d\mu_\nu(\my), \qquad  \mx\in \RR_+^d, \quad t>0,
}
where 
\eqx{
	\WW_{t,\nu}(\mx,\my) = \WW_{t,\nu_c}^{\rm{cls}}(\mx_1,\my_1) \WW_{t,-\nu_e}^{\rm{exo}}(\mx_2,\my_2), \qquad \mx = (\mx_1,\mx_2), \my=(\my_1,\my_2) \in \RR_+^{d_1}\times\RR_+^{d_2}. 
}


We define the Hardy space for the operator $\BB_\nu$ exactly in the same way as in \eqref{h1-classical}, i.e.\ by means of the maximal operator associated with the semigroup of operators $\set{\WW_{t,\nu}}$,
\eqx{
H^1(\BB_\nu) = \set{f\in L^1(\RR_+^d, d\mu_\nu) \ : \ \norm{f}_{H^1(\BB_\nu)} : = \norm{\sup_{t>0} \abs{\WW_{t,\nu}f}}_{L^1(\RR_+^d, d\mu_\nu)} < \8}.
}
We will prove that elements of $H^1(\BB_\nu)$ have an atomic decomposition, where atoms are either classical atoms \eqref{classical-atoms} or some additional atoms of the form $\mu_\nu(Q)^{-1}\mathbbm{1}_Q$ for some cube $Q\subset \RR_+^d$. Such atoms $\mu_\nu(Q)^{-1}\mathbbm{1}_Q$ are called ''local atoms'' and notice that they do not satisfy cancellation condition. 
More precisely, we will show that every function $f\in$$H^1(\BB_\nu)$ belongs also to the local atomic Hardy space $H^1_{\rm{at}}(\qq, \mu_\nu)$ associated with some family of cubes $\qq$ in $\RR_+^d$. That means that $f$ can be decomposed into a sum $f=\sum_{k} \la_k a_k$, where $\sum_k |\la_k| < \8$ and $a_k$ are either classical atoms described in \eqref{classical-atoms} and supported in cubes $Q\in\qq$ or atoms of the form $a_k = \mu_\nu(Q)^{-1}\mathbbm{1}_Q$, $Q\in\qq$.

For the general Bessel operator $\BB_\nu$ let us define the family of cubes $\qq_{\BB}$ which arises as follows. Let $\dd=\set{[2^n,2^{n+1}] \ : \ n\in\ZZ}$ be the collection of all closed dyadic intervals in $\RR_+$. Consider a family $\set{\RR_+^{d_1}\times Q_1\times \ldots \times Q_{d_2}: Q_i \in \dd}$. Then tile each cylinder $\RR_+^{d_1}\times Q_1\times \ldots \times Q_{d_2}$ with countably many cubes having diameters equal to the smallest of the diameters of $Q_1,\ldots,Q_{d_2}$. The family $\qq_{\BB}$ consists of all of these smaller cubes. For a rigorous definition of $\qq_{\BB}$ see Section \ref{sec-product}.

The main result of this paper is the following.


\thm{mix-cls-exo}{Let $d_1\geq 0$, $d_2\geq 1$. Assume that $\nu_c \in (-1,\8)^{d_1}$ and $\nu_e\in(0,\8)^{d_2}$. Then $H^1(\BB_{\nu_c}^{\rm{cls}} + \BB_{-\nu_e}^{\rm{exo}})$ and $H^1_{\rm{at}}(\qq_{\BB}, \mu_\nu)$ are isomorphic as Banach spaces.}

 In the proof we will use only general properties of $\BB_{\nu_c}^{\rm{cls}}$ and $\BB_{-\nu_e}^{\rm{exo}}$, therefore some auxiliary results in Sections \ref{sec-general} and \ref{sec-product} will be formulated in a more general context. Moreover, we introduce universal conditions on a semigroup kernel so that the product case can be deduced from a lower-dimensional information. The methods we use have roots in \cite{KaniaPlewaPreisner}, however in that paper the Lebesgue measure case is considered, and here we need to adapt them to our situation, which requires some effort. 

\subsection{Organization of the paper} Section \ref{preliminaries} contains definitions and notation used in the paper. Also, some auxiliary results are proved there. In Section \ref{sec-general} we consider a general self-adjoint and nonnegative operator $L$, which generates a semigroup of operators possessing an integral representation. We present general assumptions on the semigroup integral kernel and some family of cubes that are sufficient to prove an atomic decomposition of local type for $H^1(L)$. Section \ref{sec-product} is devoted to a similar atomic characterization in a product situation. More precisely, we slightly generalize conditions from Section \ref{sec-general} and show that if (lower-dimensional) component operators $L_1,L_2$ satisfy them, then the operator $L=L_1+L_2$ also does. Hence, one may deduce the result for the sum of operators by checking ``one-dimensional'' conditions. 
Finally, in Section \ref{sec-exotic}, we prove Theorem \ref{mix-cls-exo}. 

Throughout the paper we use standard notation. In particular, $C$ and $c$ at each occurrence denote some positive constants independent of relevant quantities. Values of $C$ and $c$ may change from line to line. Further, we write $\a \simeq \beta$ if there exists a positive constant $C$, independent of significant quantities, such that $C^{-1} \a \leq \beta \leq C \a$.

\section{Preliminaries} \label{preliminaries}

\subsection{Notation and terminology} \label{sec-notation-terminology}

We consider the metric measure space $(\RR_+^d, |\cdot|, d\mu_{\nu})$, where $|\cdot|$ stands for the standard Euclidean metric. It is well known that if $\nu\in(-1,\8)^d$, then $\mu_\nu$ possesses the doubling property, i.e.\ there exists $C>0$ such that
\eq{\label{doubling}
\mu_\nu(B(x,2r)) \leq C \mu_\nu(B(x,r)), \qquad x\in \RR_+^d, \quad r>0,}
where  $B(x,r) = \set{y\in \RR_+^d \ : \ |x-y| < r}$.

Since both the space and the measure have product structure, it is convenient to use cubes and cuboids rather than balls. Thus denote by
\eqx{Q(z;r_1,\ldots,r_d) = \set{x\in \RR_+^d \ : \ |x_i-z_i|\leq r_i \text{ for } i=1,\ldots,d}}
the cuboid centered at $z\in\RR_+^d$ having axial radii $r_1,\ldots,r_d >0$. When $r_1=\ldots=r_d=r$ the cuboid becomes a cube which we denote briefly by $Q(z,r)$. We denote by $d_Q$ the Euclidean diameter of a~cuboid $Q$.



\defn{def_covering}{
We call a family $\qq$ of cuboids in $\RR^d$ {\it an admissible covering} if there exist $C_1, C_2 >0$ such that:
\en{
\item $\RR_+^d = \bigcup_{Q\in \qq} Q$,
\item if $Q_1, Q_2 \in \qq$ and $Q_1\neq Q_2$ , then $\mu_\nu(Q_1\cap Q_2) = 0$,
\item if $Q = Q(z;r_1,\ldots,r_d)\in \qq$, then $r_i\leq C_1 r_j$ for $i,j \in \set{1,\ldots,d}$,
\item if $Q_1, Q_2 \in \qq$ and $Q_1 \cap Q_2 \neq \emptyset$, then $C_2^{-1} d_{Q_1} \leq d_{Q_2}\leq C_2 d_{Q_1}$.}}

Observe that item {\bf 3} means that admissible cuboids are uniformly bounded deformations of cubes. In fact, we shall often use only cubes. From now on we always assume that $\qq$ is an admissible covering of $\RR_+^d$.

Given a cuboid $Q$, by $Q^*$ we denote a (slight) enlargement of $Q$. More precisely, if $Q=(z;r_1,\ldots,r_d)$, then $Q^*:= Q(z; \kappa r_1,\ldots,\kappa r_d)$, where $\kappa>1$ is a fixed constant. Let $\qq$ be a given admissible covering of $\RR_+^d$. We fix $\kappa=\kappa(\qq)$ sufficiently close to 1, so that for any $Q_1,Q_2 \in \qq$,
\eq{\label{neighbours}
Q_1^{***}\cap Q_2^{***} \neq \emptyset \quad \iff \quad  Q_1 \cap Q_2 \neq \emptyset.
}
The family $\set{Q^{***}}_{Q\in \qq}$ is a finite covering of $\RR_+^d$,
\eqx{\label{finite_covering}
\sum_{Q\in \qq} \mathbbm{1}_{Q^{***}}(x) \leq C, \qquad x\in \RR_+^d.
}
For $Q\in \qq$ denote $N(Q) = \set{\wt{Q} \in \qq \ : \  \wt{Q}^{***} \cap Q^{***} \neq \emptyset,  }$ (all neighbors of the cuboid $Q$).

For the covering $\qq$ as above consider functions $\psi_Q\in C^1(\RR_+^d)$  satisfying
\eq{\label{partition}
0 \leq \psi_Q (x) \leq \mathbbm{1}_{Q^*}(x), \quad \norm{\psi_Q'}_\8 \leq C d_Q^{-1}, \quad \sum_{Q\in \qq} \psi_Q(x) = \mathbbm{1}_{\RR_+^d}(x), \qquad x\in\RR_+^d.
}
It is straightforward to see that such a family $\set{\psi_Q}_{Q\in \qq}$ exists. We call it a {\it partition of unity} related to $\qq$.

We now define suitable atoms and a local atomic Hardy space $H^1_{\rm{at}}(\qq,\mu_\nu)$ related to $\qq$.

\defn{qq-mu-atoms}{A function $a \colon \RR_+^d \to \CC$ is a $(\qq,\mu_\nu)$-atom if:

\begin{itemize}
\item[(i)] either there is  $Q\in\qq$ and a cube $K\subset Q^{*}$, such that  $$\supp \, a \subseteq K, \ \ \norm{a}_\8 \leq \mu_\nu(K)^{-1},  \ \ \int a\, d\mu_\nu = 0;$$
\item[(ii)] or there exists $Q\in\qq$ such that $$a = \mu_\nu(Q)^{-1}\mathbbm{1}_Q.$$
\end{itemize}
The atoms  as in {\it (ii)} are called {\it local atoms}.}

Having $(\qq, \mu_\nu)$-atoms at our disposal, we define the {\it local atomic Hardy space} $H^1_{\rm{at}}(\qq, \mu_\nu)$ related to $\qq$ in the standard way. Namely, a function $f$ belongs to $H^1_{\rm{at}}(\qq, \mu_\nu)$ if it has an atomic decomposition $f = \sum_k \la_k a_k$ with $\sum_k |\la_k| <\8$ and $a_k$ being $(\qq,\mu_\nu)$-atoms. The norm in $H^1_{\rm{at}}(\qq,\mu_\nu)$ is given by
$$\norm{f}_{H^1_{\rm{at}}(\qq,\mu_\nu)} = \inf \sum_k \abs{\la_k},$$
where the infimum is taken over all possible representations of $f$ as above. Note that $H^1_{\rm{at}}(\qq, \mu_\nu)$ is a Banach space.

\subsection{Auxiliary results}

\lem{measure-ball}{ Let $d=1$ and $\nu \in \RR$. The following estimates hold.
\begin{itemize}
\item[(a)] If $ \nu > -1$, then
\eqx{\mu_\nu(B(x,r)) \simeq \left( 1 \land \frac{r}{x}  \right)(x+r)^{2\nu+2}, \qquad x,r >0.}
\item[(b)] If $ \nu = -1$, then
\eqx{
\mu_\nu(B(x,r)) \simeq \log\frac{x+r}{x-r}, \qquad x>r >0.
}
\item[(c)] If $ \nu < -1$, then
\eqx{\mu_\nu(B(x,r)) \simeq \left( 1 \land \frac{r}{x}  \right)(x-r)^{2\nu+2}, \qquad x > r >0.}
\end{itemize}
}

\pr{{\bf (a)} If $r>x$, then
\spx{
	\mu_\nu(B(x,r)) = (2\nu+2)^{-1} (x+r)^{2\nu+2}.}
When $x\geq r > x/2$, we have
 \spx{
	\mu_\nu(B(x,r)) = (2\nu+2)^{-1} \left( (x+r)^{2\nu+2} - (x-r)^{2\nu+2}\right) \simeq x^{2\nu+2} \simeq (x+r)^{2\nu+2}.}
Finally, if $r\leq x/2$, then $x+r \simeq x \simeq x-r$. Therefore, applying the Mean Value Theorem,
\spx{
	\mu_\nu(B(x,r)) = (2\nu+2)^{-1} \left( (x+r)^{2\nu+2} - (x-r)^{2\nu+2}\right) \simeq  r x^{2\nu+1} \simeq \frac{r}{x} (x+r)^{2\nu+2}.}

	The proof of {\bf (b)} is straightforward. We pass to proving {\bf (c)}.
	If $x>r>x/2$, then
	 \spx{
	\mu_\nu(B(x,r)) = (-2\nu-2)^{-1} \left( (x-r)^{2\nu+2} - (x+r)^{2\nu+2}\right) \simeq (x-r)^{2\nu+2},}
	since now $2\nu+2<0$.
	When $r\leq x/2$, we proceed similarly as in the corresponding part of the proof of {\bf (a)} above. The conclusion follows.
}
\cor{m-ball}{Let $d=1$. Then
\begin{itemize}
    \item[(a)] for each $\nu\in\RR$,
\eqx{\mu_\nu(B(x,r)) \simeq r x^{2\nu+1}, \qquad 0 < r < x/2;
}
\item[(b)] for each $\nu > -1$,
\eqx{
\mu_\nu(B(x,r)) \simeq r (x+r)^{2\nu+1}, \qquad x,r > 0.
}
\end{itemize}

}
\pr{{\bf (a)} Since $r<x/2$, we have that $x+r\simeq x -r \simeq x$ and for $\nu\neq -1$ the claim readily follows from {\it (a)} and {\it (c)} of Lemma \ref{measure-ball}. To cover $\nu=-1$, use Lemma \ref{measure-ball}{\it (b)} and then observe that by the Mean Value Theorem $\log(x+r)-\log(x-r) \simeq r/x$. Part {\bf (b)} is a simple reformulation of Lemma~\ref{measure-ball}{\it (a)}.


}

As a consequence of 
the above results we obtain the following.

\cor{d-dim-measure}{Let $d\geq 1$ and $\nu\in \RR^d$. Then
\eqx{
	\mu_\nu(B(x,r)) \simeq \prod_{j=1}^d \mu_{\nu_j}(B(x_j,r)) \simeq  r^d \prod_{j=1}^d  x_j^{2\nu_j+1}, \qquad  x\in\RR_+^d, \quad 0< r < \min(x_1,\ldots,x_d)/2.
}
Moreover, if $\nu\in(-1,\8)^d$, then
\eqx{\mu_\nu(B(x,r)) \simeq \prod_{j=1}^d \mu_{\nu_j}(B(x_j,r)) \simeq  r^d \prod_{j=1}^d  (x_j+r)^{2\nu_j+1}, \qquad   x\in\RR_+^d,\quad r>0.
}}

\lem{kernel-bessel-integrable}{Let $\qq$ be an admissible covering of $\RR_+^d$. Assume that $\nu \in(-1,\8)^d$ and $Q \in \qq$. 
Let $y\in Q^*$. Then for each $\delta \in (0, \min(1/2, \nu_1+1,\ldots,\nu_d+1))$ and each $c>0$
\eq{\label{delta-positive-in}
	\int_{Q^{**}} \sup_{t>0} t^{\delta} \mu_\nu(B(x,\sqrt{t}))^{-1} \exp\left( - \frac{|x-y|^2}{ct} \right) \, d\mu_\nu(x) \leq C d_Q^{2\delta},}
and
\eq{ \label{delta-positive-out}
	\int_{(Q^{**})^c} \sup_{t>0} t^{-\delta} \mu_\nu(B(x,\sqrt{t}))^{-1} \exp\left( - \frac{|x-y|^2}{ct} \right) \, d\mu_\nu(x) \leq C d_Q^{-2\delta},
}
where the constant $C$ is independent of $Q$ and $y$.}

\pr{ 
Let $Q=Q_1\times\ldots\times Q_d$, with $d_Q\simeq d_{Q_1} \simeq \ldots \simeq d_{Q_d}$. First we prove \eqref{delta-positive-in}. To do that we need some auxiliary estimates.

{\bf Case 1:  $\nu_i\geq -1/2$}. Applying Corollary \ref{m-ball}{\it (b)} we have
\sp{ \label{Ii-positive}
	I_i & := \int_{Q_i^{**}} \sup_{t>0} t^{\delta/d} \mu_{\nu_i}(B(x_i,\sqrt{t}))^{-1} \exp\left( - \frac{|x_i-y_i|^2}{ct} \right) \, d\mu_{\nu_i}(x_i)\\
	 & \leq C \int_{Q_i^{**}} \sup_{t>0} t^{\delta/d-1/2} (x_i+\sqrt{t})^{-2\nu_i-1} \exp\left( - \frac{|x_i-y_i|^2}{ct} \right)  x_i^{2\nu_i+1} \, dx_i \\
	 & \leq C \int_{Q_i^{**}} \sup_{t>0} t^{\delta/d-1/2} \exp\left( - \frac{|x_i-y_i|^2}{ct} \right) \, dx_i \\
	& \leq C \int_{Q_i^{**}} |x_i-y_i|^{2\delta/d-1} \, dx_i \leq C d_{Q_i}^{2\delta/d}.
}
{\bf Case 2: $\nu_i\in(-1,-1/2)$}. In this case, since $\delta/d-1-\nu_i<0$, using Corollary \ref{m-ball}{\it (b)} we obtain
\sp{\label{Ii-negative}
I_i & \leq C \int_{Q_i^{**}} \sup_{0<t\leq x_i^2} t^{\delta/d-1/2} (x_i+\sqrt{t})^{-2\nu_i-1} x_i^{2\nu_i+1} \exp\left( - \frac{|x_i-y_i|^2}{ct} \right) \, dx_i \\
	 &  \quad + C \int_{Q_i^{**}} \sup_{t> x_i^2} t^{\delta/d-1/2} (x_i+\sqrt{t})^{-2\nu_i-1} x_i^{2\nu_i+1} \exp\left( - \frac{|x_i-y_i|^2}{ct} \right) \, dx_i  \\
	 & \leq C \int_{Q_i^{**}} |x_i-y_i|^{2\delta/d-1} \, dx_i\\
	 & \quad + C \int_{Q_i^{**}\cap \set{x_i\ :\ x_i/2 \leq |x_i-y_i|}} \sup_{t> 0} t^{\delta/d-1-\nu_i} x_i^{2\nu_i+1} \exp\left( - \frac{|x_i-y_i|^2}{ct} \right) \, dx_i \\
	 & \quad  + C \int_{Q_i^{**}\cap \set{x_i \ : \ x_i/2 > |x_i-y_i|}} \sup_{t> 0} t^{\delta/d-1-\nu_i} x_i^{2\nu_i+1} \exp\left( - \frac{|x_i-y_i|^2}{ct} \right) \, dx_i \\
	 & \leq  C d_{Q_i}^{2\delta/d} + C \int_{ \set{x_i\ :\ x_i \leq c d_{Q_i}}} x_i^{2\delta/d-1} \, dx_i + C \int_{Q_i^{**}} |x_i-y_i|^{2\delta/d-1} \, dx_i \\
	 & \leq C d_{Q_i}^{2\delta/d}.}

To get \eqref{delta-positive-in} we use Corollary \ref{d-dim-measure}, \eqref{Ii-positive} and \eqref{Ii-negative} obtaining
	\spx{
	& \int_{Q^{**}} \sup_{t>0} t^{\delta} \mu_\nu(B(x,\sqrt{t}))^{-1} \exp\left( - \frac{|x-y|^2}{ct} \right) \, d\mu_\nu(x) \leq  C \prod_{i=1}^d I_i \leq C \prod_{i=1}^d d_{Q_i}^{2\delta/d} \leq C d_Q^{2\delta}.	}

	Next, we prove \eqref{delta-positive-out} by induction with respect to the dimension $d$.

	{\bf Step 1.} Suppose $d=1$. We consider two cases.

	{\bf Case 1: $\nu \geq -1/2$}. Using Corollary \ref{m-ball}{\it (b)} we see that the left hand side of \eqref{delta-positive-out} is controlled by
	  \spx{ &\int_{(Q^{**})^c} \sup_{t>0} t^{-\delta-1/2} (x+\sqrt{t})^{-2\nu-1}  \exp\left( - \frac{|x-y|^2}{ct} \right) x^{2\nu+1} \,dx \\& \leq C \int_{(Q^{**})^c} \sup_{t>0} t^{-\delta-1/2} \exp\left( - \frac{|x-y|^2}{ct} \right) \, dx \\
	 & \leq C \int_{(Q^{**})^c} |x-y|^{-2\delta-1} \, dx \leq C d_{Q}^{-2\delta}.}
	 The last inequality follows from the relation $|x-y|\geq cd_Q$, since $x\in(Q^{**})^c$ and $y\in Q^*$.

	 {\bf Case 2: $\nu \in(-1,-1/2)$}. In this case we again apply Corollary \ref{m-ball}{\it (b)} to bound the left hand side of \eqref{delta-positive-out} by 

	 \spx{ & \int_{(Q^{**})^c} \sup_{t>0} t^{-\delta-1/2} (x+\sqrt{t})^{-2\nu-1}  \exp\left( - \frac{|x-y|^2}{ct} \right) \ x^{2\nu+1}\,dx \\ & \leq  C \int_{(Q^{**})^c} \sup_{0<t\leq x^2} t^{-\delta-1/2}  \exp\left( - \frac{|x-y|^2}{ct} \right) \, dx \\
	  & \quad + C \int_{(Q^{**})^c \cap\set{x \ : \ x \leq 2d_Q}} \sup_{t > x^2} t^{-\delta-1-\nu} \exp\left( - \frac{|x-y|^2}{ct} \right) \ x^{2\nu+1}\,dx\\
	   & \quad + C  \int_{(Q^{**})^c \cap\set{x \ : \ x > 2d_Q}} \sup_{t > x^2} t^{-\delta-1-\nu} \  x^{2\nu+1} \,dx\\
	   & \leq  C  \int_{(Q^{**})^c} |x-y|^{-2\delta-1}  \, dx + C  d_Q^{-2\delta-2-2\nu}\int_{0}^{2d_Q} x^{2\nu+1}  \, dx + C \int_{2d_Q}^{\8} x^{-2\delta-1} \, dx\\
	   & \leq  C d_Q^{-2\delta}.}

	  {\bf Step 2.} Let $d > 1$ be fixed. Suppose $\nu=(\nu_1,\wt{\nu}) = (\nu_1,\nu_2,\ldots,\nu_d) \in (-1,\8)^d$ and $Q = Q_1\times \wt{Q}$, where $d_Q\simeq d_{Q_1}\simeq d_{\wt{Q}}$. We use the notation $\mx = (x_1,\wt{x}), \my=(y_1,\wt{y}) \in \RR_+ \times \RR_+^{d-1}$. Assume that \eqref{delta-positive-out} holds for $\wt{Q}$, that is for every $0< \delta' < \min(1/2, \nu_2 +1, \ldots, \nu_d+1)$
	  \eq{ \label{ind-assumption}
	  \int_{(\wt{Q}^{**})^c} \sup_{t>0} t^{-\delta'} \mu_{\wt{\nu}}(B(\wt{x},\sqrt{t}))^{-1} \exp\left(-\frac{|\wt{x}-\wt{y}|^2}{ct} \right) \, d\mu_{\wt{\nu}}(\wt{x}) \leq C d_{\wt{Q}}^{-2\delta'}.
	  }
	  Take $\delta  \in (0, \min(1/2, \nu_1+1,\ldots,\nu_d+1))$. We split the set $(Q^{**})^c = S_1 \cup S_2 \cup S_2$, where 
	  \eqx{
	  S_1 = Q_1^{**} \times (\wt{Q}^{**})^c, \quad S_2 = (Q_1^{**})^c \times \wt{Q}^{**}, \quad S_3 = (Q_1^{**})^c \times  (\wt{Q}^{**})^c, 
	  } 
	  and write
	  \eqx{
	  \int_{(Q^{**})^c} \sup_{t>0} t^{-\delta} \mu_\nu(B(x,\sqrt{t}))^{-1} \exp\left( - \frac{|x-y|^2}{ct} \right) \, d\mu_\nu(x) = \int_{S_1} \ldots + \int_{S_1} \ldots + \int_{S_3} \ldots =: J_1 + J_2 + J_3. 
	  }

	  To estimate $J_1$ we apply Corollary \ref{d-dim-measure}{\it (b)}, \eqref{ind-assumption} and \eqref{delta-positive-in} with some small $\e>0$ such that $\delta+\e < \min(1/2,\nu_1 +1, \ldots,\nu_d+1)$. We get
	  \spx{
	  J_1 & \leq C \int_{Q_1^{**}} \sup_{t>0} t^{\e} \mu_{\nu_1}(B(x_1,\sqrt{t}))^{-1} \exp\left( - \frac{|x_1-y_1|^2}{ct} \right) \, d\mu_{\nu_1}(x_1)\\
     & \quad \times \int_{(\wt{Q}^{**})^c} \sup_{t>0} t^{-\delta-\e} \mu_{\wt{\nu}}(B(\wt{x},\sqrt{t}))^{-1} \exp\left(-\frac{|\wt{x}-\wt{y}|^2}{ct} \right) \, d\mu_{\wt{\nu}}(\wt{x})\\
     & \leq C d_{Q_1}^{2\e} d_{\wt{Q}}^{-2\delta-2\e} \leq C d_Q^{-2\delta}.}

	  We deal with $J_2$ similarly, but this time using \eqref{delta-positive-in} for the integral over the cube $\wt{Q}$ and \textbf{Step 1} for the integral over $Q_1$. This gives
	  \spx{
	  J_2 & \leq  C \int_{(Q_1^{**})^c} \sup_{t>0} t^{-\delta-\e}   \mu_{\nu_1}(B(x_1,\sqrt{t}))^{-1} \exp\left( - \frac{|x_1-y_1|^2}{ct} \right) \, d\mu_{\nu_1}(x_1)\\
     & \quad \times \int_{\wt{Q}^{**}} \sup_{t>0} t^{\e} \mu_{\wt{\nu}}(B(\wt{x},\sqrt{t}))^{-1} \exp\left(-\frac{|\wt{x}-\wt{y}|^2}{ct} \right) \, d\mu_{\wt{\nu}}(\wt{x})\\
     & \leq C d_{Q_1}^{-2\delta-2\e} d_{\wt{Q}}^{2\e} \leq C d_Q^{-2\delta}.}

	  Treating $J_3$ we apply \eqref{ind-assumption} and {\bf Step 1} and obtain
	  \spx{
	  J_3 & \leq C \int_{(Q_1^{**})^c} \sup_{t>0} t^{-\delta/2}   \mu_{\nu_1}(B(x_1,\sqrt{t}))^{-1} \exp\left( - \frac{|x_1-y_1|^2}{ct} \right) \, d\mu_{\nu_1}(x_1)\\
     & \quad \times \int_{(\wt{Q}^{**})^c} \sup_{t>0} t^{-\delta/2} \mu_{\wt{\nu}}(B(\wt{x},\sqrt{t}))^{-1} \exp\left(-\frac{|\wt{x}-\wt{y}|^2}{ct} \right) \, d\mu_{\wt{\nu}}(\wt{x})\\
     & \leq C d_{Q_1}^{-\delta} d_{\wt{Q}}^{-\delta} \leq C d_Q^{-2\delta}.}
	 This completes the proof.
}

\subsection{Local Hardy space}
In this subsection we consider $\nu\in(-1,\8)^d$. Let $\tau>0$ be fixed. We are interested in decomposing into atoms a function $f$ such that
	\eq{\label{local_max}
    \norm{\sup_{t\leq \tau^2} \abs{\BC{\nu} f}}_{L^1(\RR_+^d, d\mu_\nu)} <\8.
    }

In the classical case of the Laplace operator on $\RR^d$ equipped with Lebesgue measure, if one restricts the supremum to $0<t\leq \tau^2$ in the maximal operator, then one obtains an atomic space with the classical atoms complemented with atoms of the form $|B|^{-1}\mathbbm{1}_B$, where the ball $B$ has radius $\tau$, c.f. \cite{Goldberg_Duke}.
    It turns out that a similar phenomenon occurs in case of the classical Bessel operator. More precisely, \eqref{local_max} holds if and only if $f = \sum_k \la_k a_k$, where $\sum_k \abs{\la_k} <\8$ and $a_k$ are either the classical atoms or {\it local atoms at scale} $\tau$. The latter are atoms $a$ supported in a cube $Q$ of diameter comparable to $\tau$ such that $\norm{a}_\8 \leq \mu_\nu(Q)^{-1}$ but we do not impose cancellation condition. In other words one could say that these atoms built the space $H^1_{\rm{at}}(\qq^{\{\tau\}},\mu_\nu)$ introduced in Section~\ref{sec-notation-terminology}, where $\qq^{\{\tau\}}$ is a covering of $\RR_+^d$ by cubes with diameter $\tau$. The next proposition states a local atomic decomposition theorem that will be suitable for the proof of our main result. This proposition can be obtained by known methods based on Uchiyama \cite[Cor.\ 1']{Uchiyama}. We refer the reader also to \cite[Sec.\ 4]{DPW_JFAA}, where the authors check assumptions of Uchiyama's Theorem in the classical Bessel framework for the whole range of the parameter $\nu\in(-1,\8)^d$. For the sake of completeness, below we present a sketch of the proof.

    Let $\qq^{\set{\tau}}$ be a section of $\RR_+^d$ consisting of cuboids, which have diameter uniformly comparable to~$\tau$, i.e.\ there exists a positive constant $C$ such that $C^{-1} \tau \leq d_Q \leq C \tau$ for every $Q\in \qq^{\set{\tau}}$. 
 \prop{prop_local}{
 There exists $C>0$ independent of $\tau$ such that:
 \en{
 \item[(a)] For every classical atom $a$ supported in $K\subset Q^*$ or atom of the form $a=\mu_\nu(Q)^{-1}\mathbbm{1}_{Q}$, where $Q\in\qq^{\set{\tau}}$, 
 we have
 	$$\norm{\sup_{t\leq \tau^2}\abs{\BC{\nu} a}}_{L^1(\RR_+^d,d\mu_\nu)} \leq C.$$

 \item[(b)] If $\supp f \subseteq Q^*$
 and
 	$$\norm{\sup_{t\leq \tau^2}\abs{\BC{\nu} f}}_{L^1(Q^*,d\mu_\nu)}=M<\8,$$
then there exist a sequence $\la_k$ and $(\qq^{\set{\tau}},\mu_\nu)$-atoms $a_k$, such that $f = \sum_k \la_k a_k$, $\sum_k \abs{\la_k} \leq C M$, and $a_k$ are either the classical atoms supported in $Q^{*}$ or $a_k = \mu_\nu(Q)^{-1}\mathbbm{1}_{Q}$.
 }
}
\pr{ {\bf (a)} If $a$ is a classical atom satisfying cancellation condition, then the statement follows from \cite[Prop.\ 4.1]{DPW_JFAA}. So assume that $a=\mu_\nu(Q)^{-1}\mathbbm{1}_Q$. 
The maximal operator associated with $\BC{\nu}$ is bounded on $L^p(\RR_+^d,d\mu_\nu)$ for $p>1$, by Stein's general maximal theorem for semigroups of operators \cite[p.\ 73]{Stein_Topics}. Hence, using this fact and the Schwarz inequality,
\eqx{
	\norm{\sup_{t\leq \tau^2}\abs{\BC{\nu} a}}_{L^1(Q^{**},d\mu_\nu)} \leq C \mu_\nu(Q)^{1/2} \norm{\sup_{t>0}\abs{\BC{\nu} a}}_{L^2(\RR_+^d,d\mu_\nu)} \leq C \mu_\nu(Q)^{1/2} \norm{a}_{L^2(\RR_+^d,d\mu_\nu)} \leq C.
}
To estimate the norm on $(Q^{**})^c$ observe that $|x-y| \gtrsim \tau$ when $y\in Q$ and $x\in(Q^{**})^c$. Applying \eqref{classical-gauss} and \eqref{delta-positive-out} with sufficiently small $\e>0$ we obtain
\spx{
	\norm{\sup_{t\leq \tau^2} \abs{\BC{\nu} a}}_{L^1((Q^{**})^c,d\mu_\nu)} & \leq C \int_{(Q^{**})^c} \sup_{t\leq \tau^2} \int_Q \mu_\nu(Q)^{-1} \mu_\nu(B(x,\sqrt{t}))^{-1} e^{-\frac{|x-y|^2}{ct} } \, d\mu_\nu(y) d\mu_\nu(x) \\
	& \leq C \mu_\nu(Q)^{-1} \tau^{2\e} \int_Q \int_{(Q^{**})^c} \sup_{t\leq \tau^2} t^{-\e} \mu_\nu(B(x,\sqrt{t}))^{-1} e^{ -\frac{|x-y|^2}{ct} } \, d\mu_\nu(x) d\mu_\nu(y) \\
	& \leq C \mu_\nu(Q)^{-1}  \tau^{2\e} d_Q^{-2\e} \int_Q \, d\mu_\nu(y) \leq C,
	}
	since $d_Q \simeq \tau$.

	{\bf (b)} Consider $X=Q^*$ as a space. 
	Define $\la_0 = \int f d\mu_\nu$ and $g = f - \la_0 \mu_\nu(Q)^{-1}\mathbbm{1}_Q$. Notice that 
	\eq{\label{la0-norm}
	|\la_0| \leq \norm{f}_{L^1(Q^*,d\mu_\nu)} \leq M,}
	and
	\eq{\label{integral-zero}\int g \ d\mu_\nu = 0.}
	Therefore $\norm{g}_{L^1(Q^{*},d\mu_\nu)} \leq 2 \norm{f}_{L^1(Q^*,d\mu_\nu)} \leq 2M$. Moreover, $\norm{\sup_{t\leq \tau^2} \abs{\BC{\nu} g}}_{L^1(Q^*,d\mu_\nu)} \leq C M$. Indeed, using {\it (a)} and \eqref{la0-norm} we have
	\eqx{
	\norm{\sup_{t\leq \tau^2} \abs{\BC{\nu} g}}_{L^1(Q^*,d\mu_\nu)} \leq \norm{\sup_{t\leq \tau^2} \abs{\BC{\nu} f}}_{L^1(Q^*,d\mu_\nu)} + |\la_0| \norm{\sup_{t\leq \tau^2} \abs{\BC{\nu} \left( \mu_\nu(Q)^{-1}\mathbbm{1}_Q\right)}}_{L^1(Q^*,d\mu_\nu)} \leq CM.
	}
	For $t>\tau^2$, by \eqref{classical-gauss} and \eqref{doubling} we obtain
	\spx{
	\norm{\sup_{t >\tau^2} \abs{\BC{\nu} g}}_{L^1(Q^*,d\mu_\nu)} & \leq C \int_{Q^*} \sup_{t > \tau^2} \int_{Q^*} \mu_\nu(B(x,\sqrt{t}))^{-1} |g(y)| \, d\mu_\nu(y) d\mu_\nu(x) \\
	& \leq C \int_{Q^*} |g(y)| \int_{Q^*}  \sup_{t > \tau^2} \mu_\nu(B(x,\sqrt{t}))^{-1}  \, d\mu_\nu(x) d\mu_\nu(y) \\
	& \leq C \int_{Q^*} |g(y)| \int_{Q^*} \mu_\nu(B(z,\tau))^{-1}  \, d\mu_\nu(x) d\mu_\nu(y) \\
	& \leq C \norm{g}_{L^1(Q^*,d\mu_\nu)}  \\
	& \leq C M.
	}

	 Repeating the proof of \cite[Prop.\ 4.1]{DPW_JFAA}, where the authors use Uchiyama's Theorem \cite[Cor.\ 1']{Uchiyama}, we conclude that the function $g$ may be decomposed as
	\eqx{ g = \wt{\la} \mu_\nu(Q^*)^{-1} \mathbbm{1}_{Q^*} + \sum_{k} \la_k a_k, 
	} 
	where $a_k$ are atoms supported in $X=Q^*$ and satisfy cancellation condition. The constant atom $\mu_\nu(Q^*)^{-1}\mathbbm{1}_{Q^*}$ appears since $\mu_\nu(X) = \mu_\nu(Q^*) < \8$, see \cite[Sec.\ 2]{Uchiyama}. However, by \eqref{integral-zero} we get $\wt{\la} = 0$. Further, $\inf \sum_k |\la_k| \simeq \norm{\sup_{t >0} \abs{\BC{\nu} g}}_{L^1(Q^*,d\mu_\nu)} \leq CM$, which combined with \eqref{la0-norm} gives the desired conclusion.}

\section{Local atomic decomposition}\label{sec-general}

\subsection{Local atomic decomposition theorem}\label{loc_at_dec_thm}
Let $L$ be a self-adjoint and nonnegative operator defined on $L^2(\RR_+^d, d\mu_\nu)$. Denote by $T_t = \exp\left( -tL \right)$ the semigroup generated by $L$ and suppose there exists an integral kernel $T_t(x,y)$, such that $T_tf(x) = \int_{\RR_+^d} T_t(x,y) f(y) d\mu_\nu(y)$, $t>0$. Similarly as in \eqref{h1-classical} we consider the maximal Hardy space for the operator $L$
\eqx{\label{h1-L} 
H^1(L) = \set{f\in L^1(\RR_+^d, d\mu_\nu) \ : \ \norm{f}_{H^1(L)} : = \norm{\sup_{t>0} \abs{T_t f}}_{L^1(\RR_+^d, d\mu_\nu)} <\8}.
}

For $\nu\in(-1,\8)^d$, an admissible covering $\qq$, and the semigroup $\set{T_t}$ we consider the following conditions.

\eq{\label{a0} \tag{$A_0$}\begin{split}
& 0 \leq T_t(x,y) \leq C   \mu_{\nu}(B(x,\sqrt{t}))^{-1} \ \exp\left(-\frac{|x-y|^2}{ct} \right), \\
&  \qquad t>0, \quad Q\in\qq, \quad x\in N(Q), \quad y\in Q^*.
\end{split}}
\eq{\label{a1-prim} \tag{$A_1'$}
\sup_{y\in Q^{*}} \int_{(Q^{**})^c} \sup_{t>0} T_t(x,y) \, d\mu_{\nu}(x) \leq C, \qquad Q\in\qq.
}
\eq{\label{a2-prim} \tag{$A_2'$}
\sup_{y\in Q^{*}} \int_{Q^{**}} \sup_{t\leq d_Q^2} \abs{T_t(x,y) - \BC{\nu}(x,y)} \, d\mu_{\nu}(x) \leq C, \qquad Q\in\qq.}
Let us emphasize that the constants $C,c$ in these above conditions are independent of $Q\in\qq$.

The main result of this section is the following.
\thm{first-main}{Let $\nu \in (-1,\8)^d$. Assume that $L, \set{T_t}$ with an admissible covering $\qq$ satisfy \eqref{a0}, \eqref{a1-prim} and \eqref{a2-prim}. Then $H^1(L)$ and $H^1_{\rm{at}}(\qq, \mu_\nu)$ are isomorphic as Banach spaces.}

To prove Theorem \ref{first-main} we need one more lemma.


\lem{a3a4}{ Let $\nu\in(-1,\8)^d$. Assume that $\set{T_t}$ and $\qq$ satisfy conditions \eqref{a0}, \eqref{a1-prim} and a family $\set{\psi_Q}_{Q\in \qq}$ satisfies \eqref{partition}. Then 
\eq{ \label{a3}
\sup_{y \in Q^{*}} \int_{Q^{**}} \sup_{t > d_Q^2} T_t(x,y) \, d\mu_{\nu}(x) \leq C, \qquad Q\in \qq,}
\eq{ \label{a4}
\sup_{y\in \RR_+^d} \sum_{Q\in\qq} \int_{Q^{**}}  \sup_{t\leq d_Q^2} T_t(x,y) \abs{\psi_Q(x) - \psi_Q(y)} \, d\mu_{\nu}(x) < \8.}
}

\pr{
First, we prove \eqref{a3}. Let $Q\in\qq$ and take $y\in Q^*$. Notice that if $x\in Q^{**}$ then in particular $x\in N(Q)$, hence we may use \eqref{a0}. We have
\spx{
\int_{Q^{**}} \sup_{t > d_Q^2} T_t(x,y) d\mu_{\nu}(x) & \leq C \int_{Q^{**}} \sup_{t > d_Q^2} \mu_{\nu}(B(x,\sqrt{t}))^{-1} \, d\mu_{\nu}(x)\\
& \leq C  \int_{Q^{**}} \mu_{\nu}(B(x,d_Q))^{-1} \, d\mu_{\nu}(x)\\
& \leq C \mu_{\nu}(Q)^{-1}  \int_{Q^{**}} \, d\mu_{\nu}(x) \leq C,
}
since $\mu_\nu(B(x,d_Q)) \simeq \mu_\nu(Q)$, $x\in Q^*$, by the doubling property of $\mu_\nu$.

Next we show \eqref{a4}. Let $y\in \RR_+^d$. There exists $Q_0\in \qq$ such that $y\in Q_0$. Recall that $N(Q_0) = \set{Q\in\qq \ : \ Q^{***}\cap Q_0^{***} \neq \emptyset}$. We write
\spx{
\sum_{Q\in\qq} \int_{Q^{**}}  \sup_{t\leq d_Q^2} T_t(x,y) \abs{\psi_Q(x) - \psi_Q(y)} \, d\mu_{\nu}(x) = \sum_{Q\in N(Q_0)} \ldots + \sum_{Q\in\qq\setminus N(Q_0)} \ldots =: S_1 + S_2.
}

To estimate $S_2$ we use $0\leq \psi_Q \leq 1$ and \eqref{a1-prim}, getting
\spx{
S_2 & \leq \sum_{Q \in \qq\setminus N(Q_0)} \int_{Q^{**}}\sup_{t>0}  T_t(x,y) \, d\mu_{\nu}(x)\\
&  \leq C \int_{(Q_0^{**})^c}  \sup_{t>0}  T_t(x,y) \, d\mu_{\nu}(x) \leq C.}

To treat $S_1$,  write $Q=Q_1\times\ldots\times Q_d$. We claim that, given $j = 1,\ldots,d$, and $\e_j\in(0,2\nu_j+2)$, for $x_j\in N(Q_j)$, $y_j\in Q_j^*$ one has the bound
\eq{\label{sup-t}
P := \sup_{t\leq d_Q^2}  \mu_{\nu_j}(B(x_j,\sqrt{t}))^{-1} \exp\left( -\frac{|x-y|^2}{c_j t} \right) |x-y|^{\e_j} x_j^{2\nu_j+1} \leq C \min(x_j, |x_j-y_j|)^{-1+\e_j},
}
that holds with $C$ independent of $Q\in\qq$, $x_j$ and $y_j$. This will be verified in a moment. Now choose $\e_j$, $j=1,\ldots,d$, such that $\e_j\in(0,2\nu_j+2)$ and $\e:=\sum_j \e_j < 1$. Using \eqref{a0}, the Mean Value Theorem for $\psi_Q$, \eqref{partition}, Corollary \ref{d-dim-measure} and then \eqref{sup-t} we obtain 
\spx{
S_1 & \leq C \sum_{Q\in N(Q_0)} \int_{Q^{**}}  \sup_{t\leq d_Q^2} \mu_\nu(B(x,\sqrt{t}))^{-1} \exp\left( - \frac{|x-y|^2}{ct} \right) \frac{|x-y|}{d_Q} \, d\mu_{\nu}(x)\\
& \leq C \sum_{Q\in N(Q_0)} d_Q^{-\e} \int_{Q^{**}} \prod_{j=1}^d \left(  \sup_{t\leq d_Q^2} \mu_{\nu_j}(B(x_j,\sqrt{t}))^{-1} \exp\left( - \frac{|x-y|^2}{c d t} \right) |x-y|^{\e_j} x_j^{2\nu_j+1} \right) \, dx\\
& \leq C \sum_{Q\in N(Q_0)} d_Q^{-\e}   \prod_{j=1}^d \int_{Q_j^{**}}  \min(x_j,|x_j-y_j|)^{-1+\e_j} \, dx_j\\
& \leq C \sum_{Q\in N(Q_0)} d_Q^{-\e}   \prod_{j=1}^d \Bigg( \int_{\set{x_j \ : \ x_j\leq |x_j-y_j| \leq cd_Q}}  x_j^{-1+\e_j} \, dx_j + \int_{\set{x_j \ : \ x_j >  |x_j-y_j|}\cap Q_j^{**}}  |x_j-y_j|^{-1+\e_j} \, dx_j \Bigg)\\
 & \leq C \sum_{Q\in N(Q_0)} d_Q^{-\e}  \prod_{j=1}^d d_{Q_j}^{\e_j}  \leq C,}
since $d_{Q_j} \simeq d_Q$ for $j=1,\ldots,d$.

It remains to show \eqref{sup-t}. We will consider several cases and subcases.

{\bf Case 1:} $\nu_j \geq -1/2$. In this case, by Corollary \ref{m-ball}{\it (b)} we have
\spx{	
	P & \leq C \sup_{t\leq d_Q^2}  t^{-1/2} (x_j+\sqrt{t})^{-2\nu_j-1} \exp\left( -\frac{|x-y|^2}{c_j t} \right)|x-y|^{\e_j} x_j^{2\nu_j+1} =: P_1.
}

{\bf Case 1a:} $x_j < |x_j-y_j|/2$. We have 
\spx{   P_1 & \leq C \sup_{t\leq d_Q^2}  t^{-1-\nu_j} \exp\left( -\frac{|x-y|^2}{c_j t} \right) |x-y|^{\e_j} x_j^{2\nu_j+1} \\
	& \leq C   x_j^{2\nu_j+1}|x-y|^{-2\nu_j-2+\e_j} \\
	& \leq C   x_j^{2\nu_j+1}|x_j-y_j|^{-2\nu_j-2+\e_j} \\
	& \leq C  x_j^{-1+\e_j}.
} 
 
{\bf Case 1b:} $x_j \geq |x_j-y_j|/2$. Here we obtain
\spx{	P_1 & \leq C \sup_{t\leq d_Q^2}  t^{-1/2} \exp\left( -\frac{|x-y|^2}{c_j t} \right) |x-y|^{\e_j} \\
	& \leq C  |x-y|^{-1+\e_j} \\
	& \leq C  |x_j-y_j|^{-1+\e_j}.} 

{\bf Case 2:} $\nu_j \in(-1,-1/2)$. Using Corollary \ref{m-ball}{\it (b)} we get
\spx{	P  &\leq  C \sup_{t\leq x_j^2}  t^{-1/2} (x_j+\sqrt{t})^{-2\nu_j-1} \exp\left( -\frac{|x-y|^2}{c_j t} \right) |x-y|^{\e_j} x_j^{2\nu_j+1}\\
	& \quad + C \sup_{ x_j^2< t \leq d_Q^2 }  t^{-1/2} (x_j+\sqrt{t})^{-2\nu_j-1}  \exp\left( -\frac{|x-y|^2}{c_j t} \right) |x-y|^{\e_j} x_j^{2\nu_j+1}  =: P_2.} 

{\bf Case 2a:} $x_j < |x_j-y_j|/2$. We write
\spx{ P_2 &\leq  C \sup_{t > 0}  t^{-1/2} \exp\left( -\frac{|x-y|^2}{c_j t} \right) |x-y|^{\e_j} \\
	& \quad + C \sup_{x_j^2< t \leq d_Q^2 } t^{-\nu_i-1} x_j^{2\nu_j+1}  \exp\left( -\frac{|x-y|^2}{c_j t} \right) |x-y|^{\e_j}\\
 &\leq  C |x-y|^{-1+\e_j} + C x_j^{2\nu_j+1} |x-y|^{-2\nu_j-2+\e_j}\\
	&\leq C |x_j-y_j|^{-1+\e_j} + C x_j^{-1+\e_j}\\
	&\leq  C  x_j^{-1+\e_j}.} 

{\bf Case 2b:} $|x_j-y_j|/2 < x_j \leq d_Q$. In this case
\spx{P_2 &\leq C \sup_{t > 0}  t^{-1/2} \exp\left( -\frac{|x-y|^2}{c_j t} \right)  |x-y|^{\e_j}\\
	&\quad + C \sup_{x_j^2< t \leq d_Q^2 } t^{-\nu_i-1}  \exp\left( -\frac{|x-y|^2}{c_j t} \right) |x-y|^{\e_j} |x_j-y_j|^{2\nu_j+1}\\
 &\leq C |x-y|^{-1+\e_j} + C  |x-y|^{-2\nu_j-2+\e_j} |x_j-y_j|^{2\nu_j+1}\\
	&\leq  C |x_j-y_j|^{-1+\e_j}.} 

{\bf Case 2c:} $d_Q < x_j$. Here we have
\spx{ P & \leq  C \sup_{t\leq d_Q^2}  t^{-1/2} (x_j+\sqrt{t})^{-2\nu_j-1}   \exp\left( -\frac{|x-y|^2}{c_j t} \right)  |x-y|^{\e_j} x_j^{2\nu_j+1}\\
	& \leq  C \sup_{t\leq d_Q^2}  t^{-1/2} \exp\left( -\frac{|x-y|^2}{c_j t} \right)  |x-y|^{\e_j}\\\
 &\leq  C |x-y|^{-1+\e_j} \leq  C |x_j-y_j|^{-1+\e_j}.} 
 This finishes proving \eqref{sup-t}.}


\subsection{Proof of Theorem \ref{first-main} }
We shall prove the two inclusions.

\noindent{\bf First inclusion.} Let $f\in H^1_{\rm{at}}(\qq,\mu_\nu)$. We need to show $\norm{\sup_{t>0}\abs{T_t f}}_{L^1(\RR_+^d,d\mu_{\nu})} \leq C \norm{f}_{H^1_{\rm{at}}({\qq},\mu_{\nu})}$. By the standard density argument it is enough to check that $\norm{\sup_{t>0} |T_t a|}_{L^1(\RR_+^d,d\mu_{\nu})} \leq C$ for every $(\qq,\mu_{\nu})$-atom $a$, where $C$ does not depend on $a$. 

Take an atom $a$ associated with a cuboid $Q\in\qq$, see Definition \ref{qq-mu-atoms}. Using \eqref{a1-prim},  \eqref{a3}, \eqref{a2-prim} and Proposition \ref{prop_local}{\it (a)} with $\tau=d_Q$ we get

\spx{
	\norm{\sup_{t>0} |T_t a|}_{L^1(\RR_+^d,d\mu_{\nu})}  \leq &\norm{\sup_{t>0} |T_t a|}_{L^1((Q^{**})^c,d\mu_{\nu})}  + \norm{\sup_{t> d_Q^2} |T_t a|}_{L^1(Q^{**},d\mu_{\nu})} \\
    &+ \norm{\sup_{t \leq d_Q^2} |(T_t - \BC{\nu} )a|}_{L^1(Q^{**},d\mu_{\nu})}  + \norm{\sup_{t \leq d_Q^2} |\BC{\nu} a|}_{L^1(Q^{**},d\mu_{\nu})} \leq C.}

\noindent{\bf Second inclusion.} Let $f\in H^1(L)$. We shall prove $\norm{f}_{H^1_{\rm{at}}({\qq},\mu_{\nu})} \leq C \norm{\sup_{t>0}\abs{T_t f}}_{L^1(\RR_+^d,d\mu_{\nu})}$. Let $\psi_Q$ be  a partition of unity related to $\qq$,  see \eqref{partition}. Then $f = \sum_{Q \in \qq} \psi_Q f$. Denote $f_Q = \psi_Q f$  and notice that $\supp \, f_Q \subset Q^{*}$. We have
\eq{\label{ht-dec}
	\BC{\nu}  f_Q = (\BC{\nu}  - T_t)f_Q + \left( T_t f_Q - \psi_Q T_t f \right) + \psi_Q  T_t f.
}
Clearly,
\eq{ \label{psi-tt-fq}
	\sum_{Q\in\qq} \norm{ \sup_{t \leq d_Q^2} \abs{\psi_Q T_t f}}_{L^1(Q^{**},d\mu_{\nu})} \leq C \norm{\sup_{t>0} |T_t f|}_{L^1(\RR_+^d, d\mu_{\nu})}.
}
Using \eqref{a2-prim},
\eq{
	\sum_{Q\in\qq} \norm{ \sup_{t \leq d_Q^2} \abs{( \BC{\nu}  - T_t) f_Q}}_{L^1(Q^{**},d\mu_{\nu})} \leq C \sum_{Q\in \qq} \norm{f_Q}_{L^1(\RR_+^d, d\mu_{\nu})} \leq C \norm{f}_{L^1(\RR_+^d,d\mu_{\nu})}.}
By \eqref{a4},
\sp{
\label{by-a4}
	& \sum_{Q\in\qq} \norm{ \sup_{t \leq d_Q^2} \abs{T_t f_Q - \psi_Q T_t f }}_{L^1(Q^{**},d\mu_{\nu})}\\
	&\leq \sum_{Q \in \qq} \int_{\RR_+^d} \abs{f(y)} \int_{Q^{**}} \sup_{t\leq d_Q^2} T_t(x,y) \abs{\psi_Q(y) - \psi_Q(x)} \, d\mu_{\nu}(x) d\mu_{\nu}(y)\\
    &\leq C \norm{f}_{L^1(\RR_+^d,d\mu_{\nu})}.
}
Combining  \eqref{psi-tt-fq} -- \eqref{by-a4} with \eqref{ht-dec} we arrive at
\eqx{
	\sum_{Q\in\qq} \norm{ \sup_{t \leq d_Q^2} \abs{ \BC{\nu}  f_Q}}_{L^1(Q^{**}, d\mu_{\nu})} \leq C\norm{\sup_{t>0} \abs{ T_t f}}_{L^1(\RR_+^d, d\mu_{\nu})}.}
Here we have used the fact that $\norm{f}_{L^1(\RR_+^d, d\mu_\nu)} \leq C \norm{\sup_{t>0}\abs{T_t f}}_{L^1(\RR_+^d, d\mu_\nu)}$. 

Now we apply Proposition~\ref{prop_local}{\it (b)} for each $f_Q$ separately with $\tau=d_Q$. For each $f_Q$ we obtain a sequence $\la_{k}^Q$ and atoms $a_{k}^Q$ such that 
\eqx{ f = \sum_Q f_Q = \sum_{Q,k} \la_{k}^Q a_{k}^Q}
and
\eqx{\sum_Q \sum_k |\la_{k}^Q| \leq C \sum_{Q\in\qq} \norm{ \sup_{t \leq d_Q^2} \abs{ \BC{\nu}  f_Q}}_{L^1(Q^{*},d\mu_{\nu})} \leq C \norm{\sup_{t>0} \abs{T_t f}}_{L^1(\RR_+^d,d\mu_{\nu})}.}
Finally, we observe that for each $Q$ the atoms $a_{k}^Q$ obtained by Proposition \ref{prop_local} are either local atoms of the form $\mu_\nu(Q)^{-1}\mathbbm{1}_Q$ or classical atoms supported in $Q^*$, therefore they are indeed $(\qq,\mu_{\nu})$-atoms.

\section{Local atomic decomposition in the product case}\label{sec-product}

\subsection{  Product of local atomic Hardy spaces.} \label{sec-product-1}
In this section we consider operators of the form $L=L_1 + \ldots + L_N$, where each $L_i$ acts nontrivially only on the variable $x_i \in X_i = \RR_+^{d_i}$. We introduce 
conditions on the kernels of the semigroups generated by $L_i$, and admissible coverings $\qq_i$ of $X_i$, which are sufficient to prove the local atomic decomposition theorem for the Hardy space $H^1(L)$.

More precisely, we consider $X = X_1 \times \ldots \times X_N  = \RR_+^{d_1} \times \ldots \times \RR_+^{d_N} = \RR_+^d$.  We equip the space $X$ with the Euclidean metric and the measure $d\mu_\nu(x) = x^{2\nu+1} dx$, where $\nu = (\nu_1, \ldots , \nu_N)$ and $\nu_i = (\nu_{i,1}, \ldots ,\nu_{i,d_i})\in (-1,\8)^{d_i}$ for $i=1, \ldots ,N$. Assume that $L_i$ is an operator on $L^2(X_i, d\mu_{\nu_i})$, as in Section~\ref{loc_at_dec_thm}. Slightly abusing the notation we keep the symbol 
\eqx{L_i = \underbrace{I\otimes \ldots \otimes I}_{i-1 \text{ times}} \otimes L_i \otimes \underbrace{I \otimes \ldots \otimes I}_{N-i \text{ times}}} 
for the operator on $L^2(X, d\mu_\nu)$, where $I$ denotes the identity operator on the corresponding subspace, and we define
\eq{\label{Lsum}
Lf(x) = L_1f(x) + \ldots + L_Nf(x), \quad x = (x_1,\ldots,x_N) \in X.
}
We denote by $\tti{i}(x_i,y_i)$, $x_i,y_i \in X_i$, the integral kernel of $\tti{i} = \exp\eee{-tL_i}$.

Given two admissible coverings $\qq_1$ and $\qq_2$ of $X_1$ and $X_2$, respectively, we would like to produce an admissible covering on $X_1\times X_2$. However, the products $\set{Q_1 \times Q_2 \ : \ Q_1 \in \qq_1, Q_2 \in \qq_2}$ do not form an admissible covering (in general item {\bf 3} of Definition \ref{def_covering} fails). Therefore, we need a finer construction. We split each $Q=Q_1\times Q_2$ (without loss of generality let us assume that $d_{Q_1} > d_{Q_2}$) into cuboids $Q^{j}$, $j=1,\ldots,M$, such that all of them have diameters comparable to $d_{Q_2}$. 
Then the cuboids $Q^{j} = Q_1^{j}\times Q_2$, $j=1,\ldots,M,$ satisfy:
\ite{
\item $Q = \bigcup_{j=1}^M Q^{j}$,
\item $\mu_\nu(Q^{i} \cap Q^{j}) = 0$ for $i,j \in \set{1,\ldots,M}$, $i\neq j$,
\item each $Q^{j}$ satisfies condition {\bf 3} of Definition \ref{def_covering}.}

\defn{def_box_prod}{The admissible covering of $X_1\times X_2$ described in the above construction will be called the product admissible covering and denoted by $\qq_1 \boxtimes \qq_2$. }



As an example, consider the family $\dd = \set{[2^n,2^{n+1}] \ : \ n\in \ZZ}$, which is an admissible covering of $\RR_+$. The product admissible covering $\dd\boxtimes\dd$ of $\RR_+^2$ is illustrated by Figure \ref{dyadic-covering} below. 
\begin{figure}[H]
\centering
\includegraphics[scale=0.4]{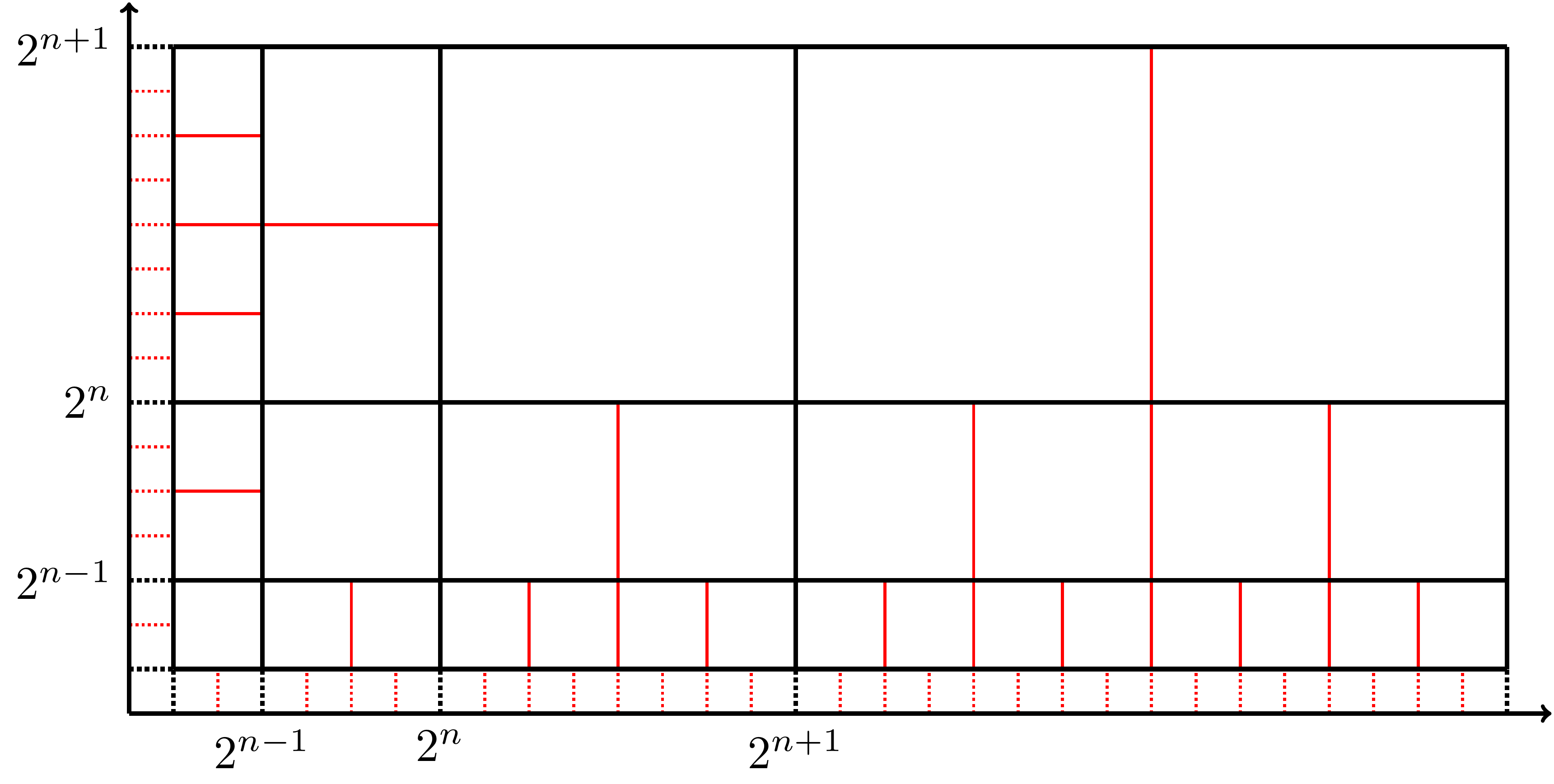}
\caption{Product admissible covering $\dd\boxtimes\dd$ of $\RR_+^2$.}\label{dyadic-covering}
\end{figure} 

For a product admissible covering $\qq_1 \boxtimes \qq_2$ we also fix a new $\kappa$ such that \eqref{neighbours} is fulfilled. The covering $\qq_1\boxtimes \ldots \boxtimes \qq_N$ of $X_1\times \ldots \times X_N$ is constructed as above by induction.

\rem{finer-split}{The collection of all the smaller cuboids $Q_1^{j}$ emerging from splitting the cuboid $Q_1$ and described before Definition \ref{def_box_prod} forms an admissible covering of $\RR_+^{d_1}$.}

 We shall assume that each $\tti{i}(x_i,y_i)$ , $i=1,\ldots,N$, and the corresponding admissible covering $\qq_i$ satisfy condition \eqref{a0}. Furthermore, we consider certain modifications of conditions \eqref{a1-prim}, \eqref{a2-prim}. Namely, for $i=1,\ldots,N$, 

 there exists $\gamma\in (0,1/3)$ such that: 
\eq{\label{a1} \tag{$A_1$}
\begin{split} \bullet \quad \text{ for each } \delta\in(-\gamma,\gamma) & \\
& \sup_{y\in Q^{*}} \int_{(Q^{**})^c} \sup_{t>0} t^{\delta} \tti{i}(x_i,y_i) \, d\mu_{\nu_i}(x_i) \leq C d_Q^{2\delta}, \qquad Q\in\qq_i, \qquad
\end{split}}
and 
\eq{\label{a2} \tag{$A_2$}
\begin{split} \bullet \quad \text{ for each } & \delta\in[0,\gamma)\\
& \sup_{y\in Q^{*}} \int_{Q^{**}} \sup_{t\leq d_Q^2} t^{-\delta} \abs{\tti{i}(x_i,y_i) - \BC{\nu_i}(x_i,y_i)} \, d\mu_{\nu_i}(x_i) \leq C d_Q^{-2\delta}, \qquad Q\in\qq_i.\end{split}}

Notice that if $\set{T_t}$ and $\qq$ satisfy \eqref{a1} and \eqref{a2}, then \eqref{a1-prim} and \eqref{a2-prim} also hold.

\rem{about-delta-negativ}{
Condition \eqref{a1} for negative $\delta\in(-\gamma,0)$ is automatically satisfied if the kernels $\tti{i}(x_i,y_i)$ satisfy the upper Gaussian estimates, that is the bound of condition \eqref{a0} for all $x_i,y_i\in X_i$. This is a direct consequence of \eqref{delta-positive-out}.
}
Our main result in this section is the following.


\thm{main-second}{ Let $L$ be as in \eqref{Lsum} and assume that for each $i=1,\ldots,N$, the kernel $\tti{i}(x_i,y_i)$ together with the corresponding admissible covering $\qq_i$ satisfy conditions \eqref{a0} -- \eqref{a2}. Then
\eqx{
H^1(L) = H^1_{\rm{at}}(\qq_1\boxtimes \ldots \boxtimes \qq_N,\mu_\nu)}
and the corresponding norms are equivalent.}

\pr{We will show the following claim. If conditions \eqref{a0} -- \eqref{a2} hold for $\tti{i}(x_i,y_i)$ together with admissible coverings $\qq_i$ for $i=1,2$, then \eqref{a0} -- \eqref{a2} also hold for $\Tt(x,y) = \tti{1}(x_1,y_1)\tti{2}(x_2,y_2)$, together with $\qq = \qq_1 \boxtimes \qq_2$. This is enough, since by simple induction we shall get that in the general case $T_t(x,y) = \tti{1}(x_1,y_1)\ldots\tti{N}(x_N,y_N)$ with $\qq_1 \boxtimes\ldots\boxtimes \qq_N$ satisfy  \eqref{a0} -- \eqref{a2}, and consequently, the assumptions of Theorem \ref{first-main} will be fulfilled.

To prove the claim let $\tti{i}(x_i,y_i)$ and $\qq_i$ satisfy \eqref{a0} -- \eqref{a2} with $\gamma_i$ for $i=1,2$. Let $0 < \gamma < \min(\gamma_1,\gamma_2, \nu_{1,1}+1,\ldots,\nu_{1,d_1}+1,\nu_{2,1}+1,\ldots,\nu_{2,d_2}+1)$. 
Suppose that $\qq \ni Q \subseteq Q_1 \times Q_2$, where $Q_1 \in \qq_1$, $Q_2 \in \qq_2$. Without loss of generality we may assume that $d_{Q_1} \geq d_{Q_2}$. Hence $Q=K\times Q_2$, where $K\subseteq Q_1$ is a cuboid emerging from the very construction of $\qq_1\boxtimes\qq_2$, see the description preceding Definition \ref{def_box_prod} and Figure \ref{prostokat} below. Denote by $z=(z_1,z_2)$ the center of $Q = K\times Q_2$. 


\begin{figure}[H]
\centering
\includegraphics[scale=0.35]{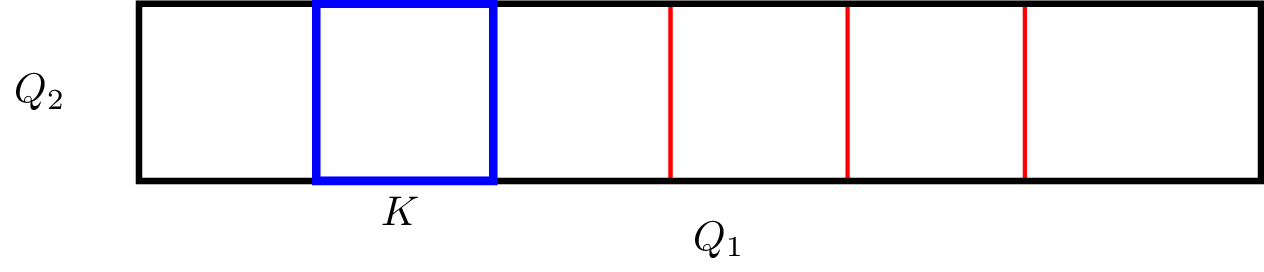}
\caption{Splitting $Q_1\times Q_2$ into admissible cuboids.}\label{prostokat}
\end{figure}

\noindent{\bf Verification of \eqref{a0}.} Notice that if $x=(x_1,x_2) \in N(Q)$, then in particular $x_1\in N(Q_1)$ and $x_2\in N(Q_2)$. Similarly, if $y = (y_1,y_2)\in Q^*$, then $y_1\in Q_1^{*}$ and $y_2\in Q_2^*$. Therefore, for such $x,y$ we may use \eqref{a0} for each of $\tti{1}(x_1,y_1)$ and $\tti{2}(x_2,y_2)$ obtaining \eqref{a0} for $T_t(x,y) = \tti{1}(x_1,y_1) \tti{2}(x_2,y_2)$.

\noindent{\bf  Verification of \eqref{a1}.} 
Fix $\delta\in(-\gamma,\gamma)$.
Let  $y\in\s{Q}{1}$. Recall that $d_Q \simeq d_K \simeq d_{Q_2} \leq d_{Q_1}$. Decompose $(Q^{**})^c = S_1\cup S_2 \cup S_3 \cup S_4 \cup S_5$, where
$$S_1 = (Q_1^{**})^c \times Q_2^{**}, \quad S_2 = K^{**} \times (Q_2^{**})^c, \quad S_3 =  \left(Q_1^{**}\setminus K^{**}\right) \times (Q_2^{**})^c,$$
$$ \quad S_4 =  (Q_1^{**})^{c} \times (Q_2^{**})^c, \quad S_5= \left( Q_1^{**}\setminus K^{**} \right) \times Q_2^{**}.$$

We start with estimating the integral over $S_1$. We shall consider two cases.

{\bf Case 1: $\delta\in(-\gamma,0]$}. In this case we take a small $\e>0$ such that $\delta-\e > -\gamma$. Using \eqref{a1} for $\tti{1}(x_1,y_1)$, \eqref{a0} for $\tti{2}(x_2,y_2)$ and then \eqref{delta-positive-in}, we have
\spx{
\int_{S_1} \sup_{t>0} t^\delta T_t(x,y) \, d\mu_{\nu}(x) &  \leq C \int_{(Q_1^{**})^c} \sup_{t>0} t^{\delta-\e} \tti{1}(x_1,y_1) \, d\mu_{\nu_1}(x_1)\\
&\quad \times \int_{Q_2^{**}} {\sup_{t>0} t^{\e} \mu_{\nu_2}(B(x_2,\sqrt{t}))^{-1}\exp\left( -\frac{|x_2-y_2|^2}{ct} \right)}\, d\mu_{\nu_2}(x_2) \\
& \leq  C d_{Q_1}^{2\delta-2\e} d_{Q_2}^{2\e} \leq C d_{Q_2}^{2\delta} \leq C d_{Q}^{2\delta},}
since $d_{Q_1} \geq d_{Q_{2}} \simeq d_Q$.

{\bf Case 2: $\delta\in(0,\gamma)$}. Then, again using \eqref{a1} for $\tti{1}(x_1,y_1)$, \eqref{a0} for $\tti{2}(x_2,y_2)$ and then \eqref{delta-positive-in} we have
\spx{
\int_{S_1} \sup_{t>0} t^\delta T_t(x,y) \, d\mu_{\nu}(x) & \leq C \int_{(Q_1^{**})^c} \sup_{t>0} \tti{1}(x_1,y_1) \, d\mu_{\nu_1}(x_1)\\
& \quad \times \int_{Q_2^{**}} \sup_{t>0} t^{\delta} \mu_{\nu_2}(B(x_2,\sqrt{t}))^{-1}\exp\left( -\frac{|x_2-y_2|^2}{ct} \right) \, d\mu_{\nu_2}(x_2) \\
& \leq  C d_{Q_2}^{2\delta} \leq C d_{Q}^{2\delta}.}

The integral over $S_2$ is treated similarly. Indeed we take a small $\e>0$ such that $\delta-\e>-\gamma$. We apply \eqref{a0} for $\tti{1}(x_1,y_1)$ and then \eqref{delta-positive-in}, and \eqref{a1} for $\tti{2}(x_2,y_2)$. We get
\spx{
\int_{S_2} \sup_{t>0} t^\delta T_t(x,y) \, d\mu_{\nu}(x) & \leq C \int_{K^{**}} \sup_{t>0}  t^{\e}\mu_{\nu_1}(B(x_1,\sqrt{t}))^{-1}\exp\left( -\frac{|x_1-y_1|^2}{ct} \right) \, d\mu_{\nu_1}(x_1)\\
&\quad \times \int_{(Q_2^{**})^c} \sup_{t>0} t^{\delta-\e}\tti{2}(x_2,y_2)\, d\mu_{\nu_2}(x_2) \\
& \leq C d_{K}^{2\e} d_{Q_2}^{2\delta-2\e} \leq C d_{Q}^{2\delta},
}
since $d_K \simeq d_{Q_2}\simeq d_Q$. By Remark \ref{finer-split} we were allowed to use Lemma \ref{kernel-bessel-integrable} for the cuboid $K$.

To estimate the integral over $S_3$ notice that since $x_1\in Q_1^{**}$ we may use \eqref{a0} for $\tti{1}(x_1,y_1)$. Let $\e>0$ be such that $\delta+\e < \gamma$. Then we use \eqref{a1} for $\tti{2}(x_2,y_2)$ and \eqref{delta-positive-out} and arrive at
\spx{ \int_{S_3} \sup_{t>0} t^{\delta} T_t(x,y) \, d\mu_{\nu}(x) & \leq  C \int_{Q_1^{**}\setminus K^{**}} \sup_{t>0} t^{-\e} \mu_{\nu_1}(B(x_1,\sqrt{t}))^{-1} \exp\left( -\frac{|x_1-y_1|^2}{ct} \right) \, d\mu_{\nu_1}(x_1)\\& \quad   \times \int_{(Q_2^{**})^c} \sup_{t>0} t^{\delta+\e} \ \tti{2}(x_2,y_2)\, d\mu_{\nu_2}(x_2) \\
 & \leq C d_K^{-2\e} d_{Q_2}^{2\delta+2\e} \leq C d_Q^{2\delta}.}

Treating the integral over $S_4$ we make use of \eqref{a1} for both $\tti{1}(x_1,y_1)$ and $\tti{2}(x_2,y_2)$. We get
\spx{\int_{S_4} \sup_{t>0} t^{\delta}T_t(x,y) \, d\mu_{\nu}(x) & \leq C \int_{(Q_1^{**})^c} \sup_{t>0} \tti{1}(x_1,y_1) \, d\mu_{\nu_1}(x_1)\\&  \quad \times \int_{(Q_2^{**})^c} \sup_{t>0} t^{\delta} \tti{2}(x_2,y_2)\, d\mu_{\nu_2}(x_2) \\
& \leq  C d_{Q_2}^{2\delta} \leq C d_Q^{2\delta}.}

Turning to the last integral over $S_5$ we consider two cases.

{\bf Case 1: $\delta\in(-\gamma,0]$}. Take $\e>0$ such that $\delta-\e > -\gamma$. We use \eqref{a0} for $\tti{1}(x_1,y_1)$ and $\tti{2}(x_2,y_2)$, and then Lemma \ref{kernel-bessel-integrable} getting
\spx{
\int_{S_5} \sup_{t>0} t^\delta T_t(x,y) \, d\mu_{\nu}(x) & \leq C \int_{Q_1^{**}\setminus K^{**}} \sup_{t>0} t^{\delta-\e} \mu_{\nu_1}(B(x_1,\sqrt{t}))^{-1}\exp\left( -\frac{|x_1-y_1|^2}{ct} \right) \, d\mu_{\nu_1}(x_1)\\
&\quad \times \int_{Q_2^{**}} {\sup_{t>0} t^{\e} \mu_{\nu_2}(B(x_2,\sqrt{t}))^{-1}\exp\left( -\frac{|x_2-y_2|^2}{ct} \right)}\, d\mu_{\nu_2}(x_2) \\ & \leq C d_{Q_2}^{2\e} \int_{ (K^{**})^c} \sup_{t>0} t^{\delta-\e} \mu_{\nu_1}(B(x_1,\sqrt{t}))^{-1}\exp\left( -\frac{|x_1-y_1|^2}{ct} \right) \, d\mu_{\nu_1}(x_1)\\
& \leq  C d_{Q_2}^{2\e} d_{K}^{2\delta-2\e}  \leq C d_{Q}^{2\delta},}
since $d_K \simeq d_{Q_{2}} \simeq d_Q$. 

{\bf Case 2: $\delta\in(0,\gamma)$}. Let $\e>0$ be such that $\delta+\e < \gamma$. Using \eqref{a0} for $\tti{1}(x_1,y_1)$, $\tti{2}(x_2,y_2)$ and then Lemma \ref{kernel-bessel-integrable} we obtain
\spx{
\int_{S_5} \sup_{t>0} t^\delta T_t(x,y) \, d\mu_{\nu}(x) & \leq C \int_{Q_1^{**}\setminus K^{**}} \sup_{t>0} t^{-\e} \mu_{\nu_1}(B(x_1,\sqrt{t}))^{-1}\exp\left( -\frac{|x_1-y_1|^2}{ct} \right) \, d\mu_{\nu_1}(x_1)\\
&\quad \times \int_{Q_2^{**}} {\sup_{t>0} t^{\delta+\e} \mu_{\nu_2}(B(x_2,\sqrt{t}))^{-1}\exp\left( -\frac{|x_2-y_2|^2}{ct} \right)}\, d\mu_{\nu_2}(x_2) \\ & \leq C d_{Q_2}^{2\delta+2\e} \int_{ (K^{**})^c} \sup_{t>0} t^{-\e} \mu_{\nu_1}(B(x_1,\sqrt{t}))^{-1}\exp\left( -\frac{|x_1-y_1|^2}{ct} \right) \, d\mu_{\nu_1}(x_1)\\
& \leq  C d_{Q_2}^{2\delta+2\e} d_{K}^{-2\e}  \leq C d_{Q}^{2\delta}.}

Condition \eqref{a1} is now verified.

\noindent{\bf Verification of \eqref{a2}.} Fix $\delta\in[0,\gamma)$.
Let $y\in\s{Q}{1}$. In this proof $\BC{\nu}$ will be the classical Bessel semigroup on $\RR_+^{d_1}$, $\RR_+^{d_2}$ or on $\RR_+^{d_1+d_2}$, depending on the context. First, notice that by \eqref{a0} and \eqref{delta-positive-in}, for any given constant $c>1$ and $i=1,2$, we have
\sp{\label{t-sim-dq}
& \int_{\s{Q_i}{2}} \sup_{c^{-1} d_{Q_i}^2\leq t \leq cd_{Q_i}^2} t^{-\gamma} \abs{T_t^{[i]}(x_i,y_i) - \BC{\nu_i}(x_i,y_i)} \, d\mu_{\nu_i}(x_i)\\
& \leq C d_{Q_i}^{-2\gamma-2\e} \int_{\s{Q_i}{2}} \sup_{c^{-1} d_{Q_i}^2\leq t \leq c d_{Q_i}^2}  t^{\e} \ \mu_{\nu_i}(B(x_i,\sqrt{t}))^{-1} \exp\left( - \frac{|x_i-y_i|^2}{c t} \right) \, d\mu_{\nu_i}(x_i)\\
& \leq C d_{Q_i}^{-2\gamma}.
 }

Since $\BC{\nu}(x,y) = \BC{\nu_1}(x_1,y_1)\BC{\nu_2}(x_2,y_2)$, by the triangle inequality we have
\spx{
& \int_{\s{Q}{2}} \sup_{t \leq d_Q^2 } t^{-\delta}\abs{\Tt(x,y) - \BC{\nu}(x,y)}d\mu_\nu(x) \leq I_1 + I_2,
}
where
\alx{
I_1 &= \int_{\s{Q}{2}} \sup_{ t \leq d_Q^2 } t^{-\delta} \tti{1}(x_1,y_1) \abs{\tti{2}(x_2,y_2) - \BC{\nu_2}(x_2,y_2)} \, d\mu_\nu(x), \\
I_2 &=  \int_{\s{Q}{2}} \sup_{ t \leq d_Q^2  } t^{-\delta} \BC{\nu_2}(x_2,y_2) \abs{\tti{1}(x_1,y_1) - \BC{\nu_1}(x_1,y_1)} \, d\mu_\nu(x).
}

Applying \eqref{a0} for $\tti{1}(x_1,y_1)$ and then \eqref{delta-positive-in}, and \eqref{a2} together with \eqref{t-sim-dq} for $\tti{2}(x_2,y_2)$, gives
\spx{
 I_1  & \leq C \int_{\s{K}{2}} \sup_{ t \leq d_Q^2  } t^{\gamma-\delta} \tti{1}(x_1,y_1) \, d\mu_{\nu_1}(x_1) \\ 
 & \quad \times \int_{Q_2^{**}} \sup_{ t \leq c d_{Q_2}^2 } t^{-\gamma} \abs{\tti{2}(x_2,y_2) - \BC{\nu_2}(x_2,y_2)} \, d\mu_{\nu_2}(x_2) \\
& \leq  C d_K^{2\gamma-2\delta} d_{Q_2}^{-2\gamma} \leq C d_Q^{-2\delta},}
since $0 \leq \delta < \gamma < \min(\gamma_1,\gamma_2)$. 

Similarly, by \eqref{classical-gauss}, \eqref{delta-positive-in}, \eqref{a2} and \eqref{t-sim-dq}, we have
\spx{
 I_2 & \leq  C \int_{\s{Q_2}{2}} \sup_{ t \leq d_Q^2  } t^{\gamma-\delta} \BC{\nu_2}(x_2,y_2) \, d\mu_{\nu_2}(x_2)\\& \quad \times \int_{\s{Q_1}{2}} \sup_{ t \leq c d_{Q_1}^2 } t^{-\gamma} \abs{\tti{1}(x_1,y_1) - \BC{\nu_1}(x_1,y_1)} \, d\mu_{\nu_1}(x_1) \\
& \leq C d_{Q_2}^{2\gamma-2\delta}d_{Q_1}^{-2\gamma} \leq C d_Q^{-2\delta},
}
since $d_{Q_1} \geq d_{Q_2} \simeq d_Q$. This finishes verification of \eqref{a2}.

The proof of Theorem \ref{main-second} is complete.}

\subsection{ Product of local and nonlocal atomic Hardy space.}

As we have mentioned, all atoms of the Hardy space $H^1(\BB_{\nu_1}^{\rm{cls}})$ satisfy cancellation condition, i.e.\ they are nonlocal atoms. However, if we consider the product $\RR_+^d = \RR_+^{d_1} \times \RR_+^{d_2}$ with the measure $\mu_\nu = \mu_{\nu_1}\otimes\mu_{\nu_2}$ and the operator $L= \BB_{\nu_1}^{\rm{cls}} + L_2$, where the semigroup $\KK_t := \exp(-tL_2)$ generated by $L_2$, together with an admissible covering $\qq_2$ satisfy conditions \eqref{a0} -- \eqref{a2} on $\RR_+^{d_2}$, then the resulting Hardy space $H^1(L)$ shall have local character.

Let $\RR_+^{d_1}\boxtimes \qq_2$ be the admissible covering that arises by splitting all the cylinders $\RR_+^{d_1} \times Q_2$, $Q_2 \in \qq_2$, into countably many cubes $Q_{1,n} \times Q_2$, where $Q_{1,n} = Q(z_n, d_{Q_2}/2)$. 
Denote by 
$\TT_t=\exp(-tL)$ the semigroup of operators generated by the operator $L$. 
There exists  an integral kernel associated with $\set{\TT_t}$ satisfying
\eqx{
\TT_t(x,y) = \BC{\nu_1}(x_1,y_1)\KK_t(x_2,y_2), \quad x=(x_1,x_2), y=(y_1,y_2) \in \RR_+^{d_1}\times\RR_+^{d_2}. }
We will prove the following.

\prop{prop-locnonloc}{
	Let $\nu=(\nu_1,\nu_2)\in(-1,\8)^{d_1}\times(-1,\8)^{d_2}$ and let $L_2$ and $L$ be as above. Assume that for $\KK_t(x_2,y_2)$ and an admissible covering $\qq_2$ of $\RR_+^{d_2}$ conditions \eqref{a0} -- \eqref{a2} hold. Then \eqref{a0} -- \eqref{a2} are satisfied for $\TT_t(x,y)$ and $\qq_L:=\RR_+^{d_1}\boxtimes \qq_2$.}

Combining Proposition \ref{prop-locnonloc} with Theorem \ref{main-second} we get the following.

\cor{cor_locnonloc}{
Let $\nu = (\nu_1,\nu_2)\in(-1,\8)^{d_1}\times(-1,\8)^{d_2}$ and let $L_2, L$ be as above. Assume that $\KK_t(x_2,y_2)$ with an
 admissible covering $\qq_2$ of $\RR_+^{d_2}$ satisfy \eqref{a0} -- \eqref{a2}. Then the spaces $H^1(L)$ and $H^1_{\rm{at}}(\RR^{d_1}\boxtimes \qq_2, \mu_\nu)$ are isomorphic as Banach spaces.}

Now we pass to proving Proposition \ref{prop-locnonloc}.

\pr{
\noindent{\bf Verification of \eqref{a0}.} Let $\qq_L \ni Q \subset\RR_+^{d_1}\times Q_2$, where $Q_2\in\qq_2$. We take $y\in Q^*$ and $x \in N(Q)$. Since $x_2\in N(Q_2)$, $y_2\in Q_2^*$, we know that \eqref{a0} holds for $\KK_t(x_2,y_2)$. Moreover, $\WW^{\rm{cls}}_{t,\nu_1}(x_1,y_1)$ satisfies \eqref{classical-gauss}, hence \eqref{a0} for $\TT_t(x,y)$ follows.

To verify \eqref{a1} and \eqref{a2} suppose that $\gamma$ is as in \eqref{a1} for $\KK_t(x_2,y_2)$ and $\qq_2$, and let \break $0< \gamma' < \min(\gamma, \nu_{1,1}+1, \ldots , \nu_{1,d_1}+1,\nu_{2,1}+1,\ldots,\nu_{2,d_2}+1)$. Let $\qq_L\ni Q = Q_1\times Q_2 \subset \RR_+^{d_1}\times Q_2$ and take $y = (y_1,y_2)\in Q^*$. Recall that $d_{Q_1}\simeq d_{Q_2} \simeq d_Q$. 

\noindent{\bf Verification of \eqref{a1}.} Let $\delta\in(-\gamma',\gamma')$. In what follows $\e>0$ will be always a fixed constant such that $\delta-\e>-\gamma'$ and $\delta+\e<\gamma'$. We decompose $(Q^{**})^c = S_1 \cup S_2\cup S_3 $, where
\eqx{
S_1 = Q_1^{**}\times(Q_2^{**})^{c}, \qquad S_2 = (Q_1^{**})^c \times Q_2^{**}, \qquad S_3 = (Q_1^{**})^c \times (Q_2^{**})^c.}

To deal with the integral over $S_1$ we consider two cases.

{\bf Case 1: $\delta\in(-\gamma',0]$}. 
We use \eqref{a1} for $\KK_t(x_2,y_2)$, \eqref{classical-gauss} for $\WW^{\rm{cls}}_{t,\nu_1}(x_1,y_1)$ and \eqref{delta-positive-in}, getting
\spx{\int_{S_1} \sup_{t>0} t^{\delta} \TT_t(x,y) \, d\mu_\nu(x) & \leq  C \int_{Q_1^{**}} \sup_{t>0} t^{\e} \mu_{\nu_1}(B(x_1,\sqrt{t}))^{-1} \exp\left( - \frac{|x_1-y_1|^2}{ct} \right) \, d\mu_{\nu_1}(x_1)\\
& \quad \times \int_{(Q_2^{**})^c} \sup_{t>0} t^{\delta-\e} \KK_t(x_2,y_2) \, d\mu_{\nu_2}(x_2)\\
& \leq  C d_{Q_1}^{2\e} d_{Q_2}^{2\delta-2\e} \leq  C d_Q^{2\delta}.}

{\bf Case 2: $\delta\in(0,\gamma')$}. By \eqref{a1} for $\KK_t(x_2,y_2)$, \eqref{classical-gauss} for $\WW^{\rm{cls}}_{t,\nu_1}(x_1,y_1)$ and \eqref{delta-positive-in} we obtain
\spx{\int_{S_1} \sup_{t>0} t^{\delta} \TT_t(x,y) \, d\mu_\nu(x) & \leq  C \int_{Q_1^{**}} \sup_{t>0} t^{\delta} \mu_{\nu_1}(B(x_1,\sqrt{t}))^{-1} \exp\left( - \frac{|x_1-y_1|^2}{ct} \right) \, d\mu_{\nu_1}(x_1)\\
& \quad \times \int_{(Q_2^{**})^c} \sup_{t>0}  \KK_t(x_2,y_2) \, d\mu_{\nu_2}(x_2)\\
& \leq C d_{Q_1}^{2\delta} \leq  C d_Q^{2\delta}.}

To estimate the integral over $S_2$ we again consider the two cases.

{\bf Case 1: $\delta\in(-\gamma,0]$}. 
Using \eqref{a0} for $\KK_t(x_2,y_2)$ and \eqref{classical-gauss} for the classical Bessel semigroup kernel, and then applying \eqref{delta-positive-out} and \eqref{delta-positive-in} we arrive at
\spx{\int_{S_2} \sup_{t>0} t^{\delta} \TT_t(x,y) \, d\mu_\nu(x) & \leq   C \int_{(Q_1^{**})^c} \sup_{t>0} t^{\delta-\e} \mu_{\nu_1}(B(x_1,\sqrt{t}))^{-1} \exp\left( - \frac{|x_1-y_1|^2}{ct} \right) \, d\mu_{\nu_1}(x_1)\\
& \quad \times \int_{Q_2^{**}} \sup_{t>0} t^{\e} \mu_{\nu_2}(B(x_2,\sqrt{t}))^{-1} \exp\left( - \frac{|x_2-y_2|^2}{ct} \right) \, d\mu_{\nu_2}(x_2)\\
& \leq  C d_{Q_1}^{2\delta-2\e} d_{Q_2}^{-2\e} \leq C d_Q^{2\delta}.}

{\bf Case 2: $\delta\in(0,\gamma')$}. 
In this case we use \eqref{a0}, \eqref{classical-gauss} and Lemma \ref{kernel-bessel-integrable}, getting
\spx{\int_{S_2} \sup_{t>0} t^{\delta} \TT_t(x,y) \, d\mu_\nu(x) & \leq   C \int_{(Q_1^{**})^c} \sup_{t>0} t^{-\e} \mu_{\nu_1}(B(x_1,\sqrt{t}))^{-1} \exp\left( - \frac{|x_1-y_1|^2}{ct} \right) \,d\mu_{\nu_1}(x_1)\\
& \quad \times \int_{Q_2^{**}} \sup_{t>0} t^{\delta+\e} \mu_{\nu_2}(B(x_2,\sqrt{t}))^{-1} \exp\left( - \frac{|x_2-y_2|^2}{ct} \right) \, d\mu_{\nu_2}(x_2)\\
 & \leq  C d_{Q_1}^{-2\e} d_{Q_2}^{2\delta+2\e} \leq C d_Q^{2\delta}.}

Finally, to treat the integral over $S_3$ we use \eqref{a1} for $\KK_t(x_2,y_2)$, \eqref{classical-gauss} for $\WW^{\rm{cls}}_{t,\nu_1}(x_1,y_1)$ and \eqref{delta-positive-out} obtaining
\spx{\int_{S_3} \sup_{t>0} t^{\delta} \TT_t(x,y) \, d\mu_\nu(x) & \leq   C \int_{(Q_1^{**})^c} \sup_{t>0} t^{-\e}\mu_{\nu_1}(B(x_1,\sqrt{t}))^{-1} \exp\left( - \frac{|x_1-y_1|^2}{ct} \right) \, d\mu_{\nu_1}(x_1)\\
& \quad \times \int_{(Q_2^{**})^c} \sup_{t>0} t^{\delta+\e} \KK_t(x_2,y_2) \, d\mu_{\nu_2}(x_2)\\
& \leq C d_{Q_1}^{-2\e} d_{Q_2}^{2\delta+2\e} \leq C d_Q^{2\delta}.}

\noindent{\bf Verification of \eqref{a2}.} Let $\delta\in[0,\gamma')$ and fix $\e>0$ such that $\delta+\e < \gamma'$. We write $\WW_{t,\nu}^{\rm{cls}}(x,y) = \WW_{t,\nu_1}^{\rm{cls}}(x_1,y_1)\WW_{t,\nu_2}^{\rm{cls}}(x_2,y_2)$. By \eqref{a2} for $\KK_t(x_2,y_2)$ and \eqref{classical-gauss} for $\WW_{t,\nu_1}^{\rm{cls}}(x_1,y_1)$ and \eqref{delta-positive-in}, we have
\spx{
& \int_{Q^{**}} \sup_{t\leq d_Q^2} t^{-\delta} \abs{\TT_t(x,y) -\WW_{t,\nu}^{\rm{cls}}(x,y) } \, d\mu_{\nu}(x)\\
& =   \int_{Q^{**}} \sup_{t\leq d_Q^2}  t^{-\delta} \ \WW_{t,\nu_1}^{\rm{cls}}(x_1,y_1) \abs{\WW_{t,\nu_2}^{\rm{cls}}(x_2,y_2) - \KK_t(x_2,y_2)} \, d\mu_{\nu}(x) \\
& \leq C \int_{Q_1^{**}}  \sup_{t\leq d_Q^2} t^{\e} \mu_{\nu_1}(B(x_1,\sqrt{t}))^{-1} \exp\left( -\frac{|x_1-y_1|^2}{ct} \right)  \, d\mu_{\nu_1}(x_1)\\
& \quad \times \int_{Q_2^{**}}  \sup_{t\leq d_Q^2} t^{-\delta-\e} \abs{\WW_{t,\nu_2}^{\rm{cls}}(x_2,y_2) - \KK_t(x_2,y_2)} \, d\mu_{\nu_2}(x_2)\\
& \leq C d_{Q_1}^{2\e} d_{Q_2}^{-2\delta-2\e} \leq C d_Q^{-2\delta}.
}

The conclusion follows.}

\section{Atomic Hardy space for the general Bessel operator}\label{sec-exotic}

In this section we prove Theorem \ref{mix-cls-exo}. Recall that we consider the operator $\BB_\nu = \BB_{\nu_c}^{\rm{cls}} + \BB_{-\nu_e}^{\rm{exo}}$ acting on functions defined on the space $X=\RR_+^{d_1}\times\RR_+^{d_2} = \RR_+^d$ equipped with the measure $\mu_{\nu} = \mu_{(\nu_c,-\nu_e)} := \mu_{\nu_c}\otimes\mu_{-\nu_e}$, where $\nu_c \in (-1,\8)^{d_1}$, $\nu_e \in (0,\8)^{d_2}$ and $d_1\geq 0$, $d_2\geq 1$. We use the notation $\nu_c = (\nu_1,\ldots,\nu_{d_1})$, $\nu_e=(\nu_{d_1+1},\ldots,\nu_{d_1+d_2})$ and $\mx = (\mx_1,\mx_2),\my=(\my_1,\my_2) \in \RR_+^{d_1}\times\RR_+^{d_2}$. Moreover, we write $\mx^{2\nu+1} = \mx_1^{2\nu_c+1}\mx_2^{-2\nu_e+1} := x_1^{2\nu_1+1} \ldots x_{d_1}^{2\nu_{d_1}+1} x_{d_1+1}^{-2\nu_{d_1+1}+1}\ldots x_{d_1+d_2}^{-2\nu_{d_1+d_2}+1}$.

\newcommand{\BEnu}{\WW_{t,-\nu}^{\rm{exo}}}

Observe that
\sp{\label{change_to_kt}
\norm{f}_{H^1(\BB_{\nu})} & = \int_X \sup_{t>0} \abs{\int_X \WW_{t,\nu_c}^{\rm{cls}}(\mx_1,\my_1)\WW_{t,-\nu_e}^{\rm{exo}} (\mx_2,\my_2) f(\my) d\mu_{(\nu_c,-\nu_e)}(\my)}  \, d\mu_{(\nu_c,-\nu_e)}(\mx)  \\
& = \int_X \sup_{t>0} \abs{\int_X \WW_{t,\nu_c}^{\rm{cls}}(\mx_1,\my_1) \mx_2^{-4\nu_e}\WW_{t,-\nu_e}^{\rm{exo}} (\mx_2,\my_2) \my_2^{-4\nu_e}f(\my)  d\mu_{(\nu_c,\nu_e)}(\my)}   d\mu_{(\nu_c,\nu_e)}(\mx)\\
& = \int_X \sup_{t>0} \abs{\int_X \WW_{t,\nu_c}^{\rm{cls}}(\mx_1,\my_1) \KK_t(\mx_2,\my_2) \wt{f}(\my) d\mu_{(\nu_c,\nu_e)}(\my)}  \, d\mu_{(\nu_c,\nu_e)}(\mx)\\
& = \int_X \sup_{t>0} \abs{\int_X\TT_t(\mx,\my) \wt{f}(\my)  d\mu_{(\nu_c,\nu_e)}(\my)}  \, d\mu_{(\nu_c,\nu_e)}(\mx)\\
& = \norm{\sup_{t>0}|\TT_t \wt{f}|}_{L^1(X,d\mu_{(\nu_c,\nu_e)})},} 
where $\wt{f}(\my) = \my_2^{-4\nu_e} f(\my)$, $\TT_t(\mx,\my) = \WW_{t,\nu_c}^{\rm{cls}}(\mx_1,\my_1) \KK_t(\mx_2,\my_2)$, and
\spx{
\KK_t(\mx_2,\my_2) & = \mx_2^{-4\nu_e} \WW_{t,-\nu_e}^{\rm{exo}} (\mx_2,\my_2) \\
& = \prod_{i=d_1+1}^{d_1+d_2} (2t)^{-1}y_i^{\nu_i} x_i^{-3\nu_i} I_{\nu_i}\left( \frac{x_iy_i}{2t} \right) \exp\left( -\frac{x_i^2+y_i^2}{4t}\right)\\
& =: \prod_{i=d_1+1}^{d_1+d_2} K_{t,\nu_i}(x_i,y_i).
}

From the above we see that $f\in H^{1}(\BB_{\nu})$ if and only if $\sup_{t>0} |\TT_t \wt{f}| \in L^1(X,d\mu_{(\nu_c,\nu_e)})$. To prove Theorem \ref{mix-cls-exo} we need the following.

\prop{a0a1a2_for_kt}{
	Let $\dd = \set{[2^n,2^{n+1}] \ : \ n\in \ZZ}$ be the admissible covering of $\RR_+$. Assume $\nu\in(0,\8)$ and let 
	\eq{ \label{kt-kernel}
	K_{t,\nu}(x,y) = (2t)^{-1}y^{\nu} x^{-3\nu} I_{\nu}\left( \frac{xy}{2t} \right) \exp\left( -\frac{x^2+y^2}{4t}\right).} Then $K_{t,\nu}$ with $\dd$ satisfy conditions \eqref{a0} -- \eqref{a2}.}

In the proof we will need the following standard asymptotics of the modified Bessel function $I_\tau$, see e.g.\ \cite[p.\ 203-204]{Watson},
\begin{align} \label{Bessel-function-small}
&&I_\tau(x) & = \Gamma\left( \tau+1\right)^{-1} \left(x/2\right)^\tau + O(x^{\tau+2}), & \text{ for } & x\sim 0,&&\\
\label{Bessel-function-large}
&&I_\tau(x) & = (2\pi x)^{-1/2} e^x + O(x^{-3/2}e^x), & \text{ for } & x \sim \8.&&
\end{align}

\pr{To verify \eqref{a0} -- \eqref{a2} for $K_{t,\nu}$ with $\dd$, let $\dd\ni Q = [2^n,2^{n+1}]$ for some $n\in \ZZ$ and take $y\in Q^*$.

\noindent{\bf Verification of \eqref{a0}.} We consider $x\in N(Q)$. Notice that then $x^{-1}y \simeq 1$, hence $K_{t,\nu}(x,y) = (x^{-1}y)^{2\nu} W_{t,\nu}^{\rm{cls}}(x,y) \simeq W_{t,\nu}^{\rm{cls}}(x,y)$. Therefore
\spx{
0 \leq K_{t,\nu}(x,y) \leq C \mu_\nu(B(x,\sqrt{t}))^{-1} \exp \left( - \frac{|x-y|^2}{ct} \right),
}
since the kernel $W_{t,\nu}^{\rm{cls}}(x,y)$ is non-negative and satisfies the Gaussian estimate \eqref{classical-gauss}.

\noindent{\bf Verification of \eqref{a1}.} Suppose $0 < \gamma < \min(1/2,\nu)$ and take $\delta\in(-\gamma,\gamma)$.
We write
\spx{\int_{(Q^{**})^c} \sup_{t>0} t^{\delta} K_{t,\nu}(x,y) \, d\mu_{\nu}(x) \leq  & \int_{(Q^{**})^c} \sup_{t>xy/2} t^{\delta} K_{t,\nu}(x,y) \, d\mu_{\nu}(x) \\
& + \int_{(Q^{**})^c\cap(0,2d_Q)} \sup_{0<t\leq xy/2} t^{\delta} K_{t,\nu}(x,y) \, d\mu_{\nu}(x)\\
& + \int_{(Q^{**})^c\cap(2d_Q,\8)} \sup_{0<t\leq xy/2} t^{\delta} K_{t,\nu}(x,y) \, d\mu_{\nu}(x)\\
 =: & I_1 + I_2 + I_3.}

To treat $I_1$ we use \eqref{kt-kernel} and \eqref{Bessel-function-small}.
\spx{I_1 & \leq C \int_0^{\8} \sup_{t>xy/2} t^{\delta} (2t)^{-1}y^{\nu} x^{-3\nu} \left( \frac{xy}{2t} \right)^{\nu} \exp\left( -\frac{x^2+y^2}{4t}\right) x^{2\nu+1} \, dx\\
& \leq C y^{2\nu} \int_0^{\8} (x^2+y^2)^{\delta-\nu-1} x \, dx\\
& = C y^{2\nu} \lim_{\e\to0} \abs{ (\e^{-2}+y^2)^{\delta-\nu} - (\e^2+y^{2})^{\delta-\nu}}\\
& = C y^{2\nu} y^{2\delta-2\nu} \leq C d_Q^{2\delta},
}
since $\delta-\nu <0$ and $y\simeq d_Q$.


To deal with $I_2$ we use \eqref{kt-kernel} and \eqref{Bessel-function-large}. Observing that $xy\leq cd_Q^2 \leq c |x-y|^2$, we get
\spx{I_2  & \leq C \int_{(Q^{**})^c\cap(0,2d_Q)} \sup_{0 < t \leq xy/2} t^{\delta-1/2} y^{\nu-1/2} x^{-\nu+1/2} \exp\left( - \frac{|x-y|^2}{4t} \right) \, dx\\
& \leq C y^{\delta+\nu-1} \int_{(Q^{**})^c\cap(0,2d_Q)} x^{\delta-\nu} \exp\left( - \frac{|x-y|^2}{cxy} \right) \, dx\\
& \leq C y^{\delta+\nu-1} \int_{(Q^{**})^c\cap(0,2d_Q)} x^{\delta-\nu} \left(\frac{|x-y|^2}{cxy} \right)^{-N} \, dx\\
& \leq C y^{\delta+\nu-1+N} d_Q^{\delta-\nu+N} \int_{(Q^{**})^c} |x-y|^{-2N}  \, dx\\
& \leq C y^{\delta+\nu-1+N} d_Q^{\delta-\nu-N+1} \leq C d_Q^{2\delta},
}
since $y\simeq d_Q$, and we choose $N$ large enough so that $\delta-\nu+N \geq 0$. 

Treating $I_3$ we observe that $x\geq 2d_Q$ implies $|x-y| \simeq x$. Then using \eqref{kt-kernel} and \eqref{Bessel-function-large} we obtain
\spx{
I_3 & \leq C \int_{(Q^{**})^c\cap(2d_Q,\8)} \sup_{0 < t \leq xy/2} t^{\delta-1/2} y^{\nu-1/2} x^{-\nu+1/2} \exp\left( - \frac{|x-y|^2}{4t} \right) \, dx\\
& \leq C y^{\nu-1/2}\int_{(Q^{**})^c\cap(2d_Q,\8)} \sup_{0 < t \leq xy/2} t^{\delta-1/2}  x^{-\nu+1/2} \left(\frac{|x-y|^2}{4t} \right)^{-N} \, dx \\
& \leq C y^{\nu-1+\delta+N}\int_{(Q^{**})^c\cap(2d_Q,\8)}  x^{\delta+N-\nu} |x-y|^{-2N} \, dx \\
& \leq C y^{\nu-1+\delta+N}\int_{(Q^{**})^c} |x-y|^{\delta-N-\nu} \, dx \\
& \leq C  y^{\nu-1+\delta+N} d_Q^{\delta-N-\nu+1} \leq C d_Q^{2\delta},
}
by taking $N$ sufficiently large and since $y\simeq d_Q$.

\noindent{\bf Verification of \eqref{a2}.} Let $0<\gamma < \min(1/2,\nu)$ and take $\delta\in[0,\gamma)$. Notice that here $x\simeq y \simeq d_Q$. Moreover, $t\leq d_Q^2 \leq c xy $. Using \eqref{kt-kernel}, \eqref{classical-kernel}, \eqref{Bessel-function-large} and the Mean Value Theorem, we obtain
\spx{& \int_{Q^{**}} \sup_{t\leq d_Q^2} t^{-\delta} \abs{K_{t,\nu}(x,y) - W_{t,\nu}^{\rm{cls}}(x,y)} \, d\mu_{\nu}(x)\\
& \leq \int_{Q^{**}} \sup_{t\leq d_Q^2} t^{-\delta}(2t)^{-1} I_{\nu}\left( \frac{xy}{2t} \right) \exp\left( -\frac{x^2+y^2}{4t}\right) x^{-\nu} y^{\nu} \abs{ x^{-2\nu} - y^{-2\nu}} x^{2\nu+1} \, dx\\
& \leq C \int_{Q^{**}} \sup_{t>0} t^{-1-\delta}\left( \frac{xy}{2t} \right)^{-1/2} \exp\left( -\frac{|x-y|^2}{4t}\right) x^{\nu+1} y^{\nu} \abs{ x - y} \xi_{x,y}^{-2\nu-1}  \, dx\\
& \leq C y^{2\nu} \int_{Q^{**}} |x-y|^{-1-2\delta} \abs{ x - y} d_Q^{-2\nu-1}  \, dx\\
& \leq C d_Q^{-1} \int_{Q^{**}} |x-y|^{-2\delta} \, dx \leq C d_Q^{-2\delta},
}
since $y\simeq d_Q$ and $\delta < 1/2$; here $\xi_{x,y}$ denotes a point between $x$ and $y$.
}

\noindent{\bf Proof of Theorem \ref{mix-cls-exo}.} Recall that $\qq_\BB = \RR_+^{d_1} \boxtimes \underbrace{\dd\boxtimes \ldots \boxtimes \dd}_{d_2 \text{ times }}$. In view of \eqref{change_to_kt} we have 
\eqx{
	\norm{f}_{H^1(\BB_{\nu})} = \norm{\sup_{t>0} |\TT_t\wt{f}|}_{L^1(X,d\mu_{(\nu_c,\nu_e)})},}
where $\wt{f}(\my) = \my_2^{-4\nu_e}f(\my)$. Combining Proposition \ref{a0a1a2_for_kt} with Proposition \ref{prop-locnonloc} and Theorem \ref{main-second} we infer that $\sup_{t>0} |\TT_t\wt{f}| \in L^1(X,d\mu_{(\nu_c,\nu_e)})$ if and only if there exist a sequence $\la_k$ and $(\qq_{\BB},\mu_{(\nu_c,\nu_e)})$-atoms $\wt{a}_k$ such that 
\eqx{\wt{f} = \sum_k \la_k \wt{a}_k \quad  \text{and} \quad \norm{\sup_{t>0} |\TT_t\wt{f}|}_{L^1(X,d\mu_{(\nu_c,\nu_e)})} \simeq \inf \sum_k |\la_k|.} 
Hence, $f$ belongs to $H^1(\BB_{\nu})$ if and only if it has a decomposition
\eqx{f(\my) = \sum_k \la_k\wt{a}_k(\my) \my_2^{4\nu_e} =: \sum_k \la_k a_k(\my).}
Thus it suffices to show that functions of the form $a_k(\my) = \wt{a}_k(\my) \my_2^{4\nu_e}$ are $(\qq_\BB,\mu_{(\nu_c,-\nu_e)})$-atoms after a~modification by multiplying by some function which is comparable to a constant. 

Suppose first that $\wt{a}_k = \mu_{(\nu_c,\nu_e)}(Q)^{-1}\mathbbm{1}_Q$ for some cube $\qq_\BB \ni Q = Q_1 \times Q_2 \subset \RR_+^{d_1} \times [2^{n_1},2^{n_1+1}]\times \ldots \times [2^{n_{d_2}},2^{n_{d_2}+1}]$. Without loss of generality we may assume that $n_1 = \min_i{n_i}$ and that $d_Q \simeq 2^{n_1} \simeq d_{Q_1} \simeq d_{Q_2} $. Denote by $z$ the center of $Q$. Then
\eqx{
a_k(\my) = \my_2^{4\nu_e} \wt{a}_k(\my) = \my_2^{4\nu_e} \frac{\mu_{-\nu_e}(Q_2)}{\mu_{\nu_e}(Q_2)} \mu_{\nu_c}(Q_1)^{-1} \mu_{-\nu_e}^{-1}(Q_2) \mathbbm{1}_Q(\my) = \my_2^{4\nu_e} \frac{\mu_{-\nu_e}(Q_2)}{\mu_{\nu_e}(Q_2)} \mu_{(\nu_c,-\nu_e)}(Q)^{-1} \mathbbm{1}_Q(\my),
} 
where $\mu_{(\nu_c,-\nu_e)}(Q)^{-1} \mathbbm{1}_Q(\my)$ is $(\qq_\BB,\mu_{(\nu_c,-\nu_e)})$-atom. Using Corollary \ref{m-ball}{\it (b)}, for $\my_2\in Q_2$ we obtain
\sp{\label{measures-cube-comparable} \my_2^{4\nu_e} \frac{\mu_{-\nu_e}(Q_2)}{\mu_{\nu_e}(Q_2)} & \simeq \my_2^{4\nu_e} \frac{ \prod_{i=d_1+1}^{d_1+d_2} \mu_{-\nu_i}( B(z_i,2^{n_1})) }{\prod_{i=d_1+1}^{d_1+d_2}\mu_{\nu_i}(B(z_i,2^{n_1}))} \simeq \prod_{i=d_1+1}^{d_1+d_2} y_i^{4\nu_i} \frac{ \prod_{i=d_1}^{d_1+d_2} z_i^{-2\nu_i+1} \cdot 2^{n_1} }{\prod_{i=d_1+1}^{d_1+d_2} z_i^{2\nu_i+1} \cdot 2^{n_1}}\\
& = \prod_{i=d_1+1}^{d_1+d_2} y_i^{4\nu} z_i^{-4\nu_i} \simeq 1,}
since $y_i \simeq z_i$ for $i=d_1+1,\ldots,d_1+d_2$. Notice that the above estimate is uniform in $Q_2$. 

On the other hand, if $\wt{a}_k$ is a $(\qq_\BB,\mu_{(\nu_c,\nu_e)})$-atom of the form $(i)$ from Definition \ref{qq-mu-atoms} associated with a~cube $\qq_\BB \ni Q = Q_1 \times Q_2 \subset \RR_+^{d_1} \times [2^{n_1},2^{n_1+1}]\times \ldots \times [2^{n_{d_2}},2^{n_{d_2}+1}]$ as above, then the function $a_k(\my) = \my_2^{4\nu_e}\wt{a}_k(\my)$ satisfies:
\en{
\item $\supp \ a_k \subseteq K = K_1\times K_2 \subseteq Q^*,$
\item similarly as in \eqref{measures-cube-comparable},
\spx{
\norm{a_k}_\8 & = \esssup_{\my\in K} \my_2^{4\nu_e} \abs{\wt{a}_k(\my)} \leq \esssup_{\my\in K} \my_2^{4\nu_e} \norm{\wt{a}_k}_\8 \\
& \leq \mu_{(\nu_c,-\nu_e)}(K)^{-1} \esssup_{\my_2\in K_2} \my_2^{4\nu_e} \frac{\mu_{-\nu_e}(K_2)}{\mu_{\nu_e}(K_2)} \leq C \mu_{(\nu_c,-\nu_e)}(K)^{-1},}
where the constant $C$ does not depend on $K_2$,
\item 
\eqx{\int a_k(\my) d\mu_{(\nu_c,-\nu_e)}(\my) = \int \wt{a}_k(\my) d\mu_{(\nu_c,\nu_e)}(\my) = 0.}
}
This finishes the proof.
\hfill \qedsymbol


\bibliographystyle{amsplain}        

\def\cprime{$'$}
\providecommand{\bysame}{\leavevmode\hbox to3em{\hrulefill}\thinspace}
\providecommand{\MR}{\relax\ifhmode\unskip\space\fi MR }
\providecommand{\MRhref}[2]{%
  \href{http://www.ams.org/mathscinet-getitem?mr=#1}{#2}
}
\providecommand{\href}[2]{#2}

\end{document}